\documentstyle[11pt,psfig]{article}
\hsize \textwidth      
\advance \hsize by -\marginparwidth
\evensidemargin \oddsidemargin  

\advance\hoffset by 5mm

\newtheorem{thm}{Theorem}
\newtheorem{prop}{Proposition}
\newtheorem{lemma}{Lemma}
\newcommand{\Rm}{I\!\!R}

\begin{document}
\title{Bulk Burning Rate in Passive - Reactive Diffusion.\\}

\author{Peter Constantin \and Alexander Kiselev \and Adam Oberman 
\and Leonid Ryzhik\\
Department of Mathematics\\  University of Chicago\\
Chicago IL 60637}


\maketitle

\begin{abstract}
We consider a passive scalar that is advected by a prescribed mean
zero divergence-free velocity field, diffuses, and reacts according
to a KPP-type nonlinear reaction. We introduce a quantity, the
bulk burning rate, that makes both mathematical and physical sense
in general situations and extends the often ill-defined notion of
front speed. We establish rigorous lower bounds for the bulk
burning rate that are linear in the amplitude of the advecting
velocity for a large class of flows. These "percolating" flows
are characterized by the presence of tubes of streamlines 
connecting distant regions of burned and unburned material and generalize shear flows. The bound contains
geometric information on the velocity streamlines and degenerates
when these oscillate on scales that are finer than the width of
the laminar burning region. We give also examples of very different
kind of flows, cellular flows with closed streamlines, and  rigorously
prove that these can produce only sub-linear enhancement of the
bulk burning rate. 
\end{abstract}

\section{Introduction}

Quite often mixtures of reactants interact in a burning region
that has a rather complicated spatial structure but is thin across. 
This reaction region moves towards the unburned reactants leaving behind the
burned ones. When the reactants are carried by an ambient fluid
then the burning rate may be enhanced. The physical reason 
for this observed speed-up is believed to be that fluid
advection tends to increase the area available for reaction.

Many important engineering applications of
combustion operate in the presence of turbulent advection, 
and therefore the influence of advection on burning has been
studied extensively by physicists, engineers and mathematicians. 
In the physics literature one can find a number of models and
approaches that yield different predictions --
relations between the turbulent intensity and the burning rate
(\cite{CW,ker-1,KA,Y}). These results are usually obtained using
heuristic models and physical reasoning. For a recent review of some
of the physics literature we refer to \cite{curved-fronts,Ro}.

The key question we wish to address is: what characteristics of the
ambient fluid flow are responsible for burning rate enhancement?
The question needs first to be made precise, because the reaction region
may be complicated and, in general, may move with an ill-defined
velocity. 

In this paper we  will define in an unambiguous fashion a quantity $V$ 
representing the bulk burning rate. $V$
makes both mathematical and physical sense in general; we study 
its relation to the advecting
velocity field in a simple model. We
provide explicit estimates of $V$ in terms of 
the magnitude of the advecting velocity and the geometry of streamlines.
We are mostly interested in the regime where the advection is 
strong but our estimates are valid for all 
values of  physical parameters, 
and do not involve any passage to limit. They are also valid for certain  
advection velocities without symmetry. 
In situations where traveling waves are
known to exist, the estimates we derive 
provide automatically bounds for the speed of the traveling waves.

The main result of this 
paper is the identification of a class of flows that are particularly 
effective 
in speeding up the bulk burning rate. 
We call these ``percolating flows'' because their main feature is the presence of
tubes of 
streamlines connecting distant regions of burned and unburned material. 
For such flows we obtain an optimal linear enhancement bound
$$
V\ge KU
$$
where $U$ represents the magnitude of the advecting velocity and
$K$ is  a proportionality factor that depends on the geometry of
streamlines but not the speed of the flow.
Other flows and in particular cellular flows, which have 
closed streamlines, on the other hand, may produce a weaker
enhancement.

We will take an analytic approach. Numerical work
concerning our results will be published elsewhere.
We consider a well-established simple model, the passive-reactive
diffusive scalar equation 
\begin{equation}
  \label{eq:1.1}
  \frac{\partial T}{\partial t} + u(x,y,t)\cdot\nabla T=\kappa\Delta T+\frac{v_0^2}{4\kappa}f(T).
\end{equation}
In this equation $0 \leq T \leq 1$ represents normalized temperature, $u$
the advecting velocity, $\kappa$ thermal 
diffusivity, $f$ the reaction and $v_0$ the laminar front speed. The advecting
velocity is divergence-free and prescribed. No feedback of $T$ on $u$ is
allowed in this simple model: $T$ is passive. The normalization is
such that the reaction rate is $\frac{v_0^2}{4\kappa}$. This is chosen
so that, in the absence of advection ($u=0$) and given (\ref{eq:1.4}) below, there
will exist reaction-diffusive laminar traveling wave
fronts that move with speed at least $v_0$ \cite{KPP}.

The equation (\ref{eq:1.1}), derived under assumptions of
approximately constant density  \cite{CW} and approximately unity
Lewis number (e.g. \cite{Ber-Lar-Roq})
 is also used 
to model 
problems in biology (\cite{Fi}), 
chemistry, and has other applications (\cite{curved-fronts}) but
certainly does not capture all the physical instabilities present in
turbulent combustion (\cite{Siv-1,Siv-2}).

The type of nonlinearity $f(T)$ we consider in this paper is
concave KPP:
\begin{equation}
  \label{eq:1.4}
f \in C^2, \,\,\,f(0)=f(1)=0, \,\,\,f'(0)=1, \,\,\,f''(x)<0. 
\end{equation}
The prototype non-linearity of the KPP type is $f(T)=T(1-T),$ called KPP 
after pioneering work of Kolmogorov, Petrovskii and Piskunov \cite{KPP}.

Other important 
types of nonlinearity $f(T)$ are the Arrhenius-type 
\[ f(T)= (1-T) e^{-\frac{A}{T}} \]
and ignition type 
\[
  f(T)= 0 {\hbox { for $T\notin(\theta,1)$}},~~
f(T)>0 \hbox { for $T\in (\theta,1)$}, ~~\theta\in(0,1). 
\]

The mathematical study of equation (\ref{eq:1.1}) 
concentrated mainly on two issues: existence of
traveling waves and asymptotic speed, and 
the homogenization regime
$\kappa\to 0$.  Traveling waves in one dimension with $u=0$ were
studied in the classical works \cite{KPP} and \cite{Fisher} with
their global asymptotic analysis addressed later in \cite{Kanel} for
the ignition nonlinearity, in \cite{AW} in higher dimensions, and in
\cite{xinx} for variable diffusivity and ignition nonlinearity.
 Traveling waves with $u \ne 0$ were shown to exist for
shear flows for KPP as well as for a more general class of nonlinearities
\cite{Ber-Lar-Lions,Ber-Nir-1,Ber-Nir-2}.  
Their stability was established in \cite{Ber-Lar-Roq,Roq-1,Mal-Roq}.  
Finally, traveling waves for periodic flows $u(x,y)$ and ignition nonlinearity as well as their
stability were studied in \cite{jxin-1,jxin-2}.
Probabilistic methods were applied for the analysis of the KPP fronts
and proof of the existence of the asymptotic speed of propagation was 
given in \cite{freidlin-4}
for periodic $u\ne 0$. To the best of our knowledge,
until now there have been no explicit estimates on the speed of propagation of
traveling waves or asymptotic speed of propagation with $u \ne 0,$
except for the perturbative small $u$ result of \cite{px}. 

The homogenization regime 
$\kappa \to 0,$ when the front width goes to zero, was  
extensively studied for KPP-type nonlinearity and 
for advection velocity that is periodic and varies either on the
integral scale \cite{freidlin-1,freidlin-2,freidlin-3} or on a small
scale that is larger or comparable to that of the front width
\cite{freidlin-4,MS,Emb-Maj-Sou-1,Emb-Maj-Sou-2,McL-Zhu}.  Recently a
similar result was established for random statistically homogeneous
ergodic advection velocities \cite{S}. A thorough
review of most of these results, both for traveling waves and
homogenization techniques is given in \cite{jxin-3}.
The result of homogenization procedures is an effective equation 
valid in  the limit 
$\kappa \to 0.$ The effective equation is typically a non-trivial
Hamilton-Jacobi equation  
\cite{Emb-Maj-Sou-1,Emb-Maj-Sou-2,MS}.
Homogenization is usually very efficient when a mean field captures
the essence of the question asked. When that is not the case
important information is lost in the limit. 


For simplicity of exposition
we will consider the reaction-diffusion-advection equation 
(\ref{eq:1.1})  in a two-dimensional strip  
\[
\Omega=\left\{x,y:~x\in(-\infty,\infty),~
  y\in[0,H]\right\},
\]
but the methods 
we introduce work in any dimension and for more general classes of domains. 
The boundary conditions are either Neumann 
\begin{equation}
  \label{eq:1.2}
  T_y(x,0)=T_y(x,H)=0,
\end{equation}
or periodic in $y$:
\begin{equation}
  \label{eq:1.2.1}
  T(x,y,t)=T(x,y+H,t).
\end{equation}
The flow $u=(u_1,u_2)$ is incompressible:
\begin{equation}
  \label{eq:1.10}
  \nabla\cdot u=0
\end{equation}
and has zero normal component at the boundary
\begin{equation}
  \label{eq:1.11}
  u\cdot n=0~~\hbox{on}~~\partial\Omega
\end{equation}
in the case of Neumann boundary conditions.
Furthermore, we assume that the total flow through the strip is zero:
\begin{equation}\label{eq:1.11.1}
\int\limits_0^H u_1(x,y)dy=0
\end{equation}
to eliminate the drift caused by the mean flow.
We assume $1 \geq T \geq 0,$ 
$T=1$ being the stable (burned) 
state of the system, while $T=0$ is unstable (unburned) state.
The equation has maximum principle \cite{MAX}, so 
if the initial data is in the $[0,1]$ range then the
solution remains in the same range for all times. 
We also assume that the solution is localized, that is
\begin{eqnarray}\label{eq:1.3}
  T(x,y,t)=1-O(e^{\lambda x})~~\hbox{for}~~x<0,~~
  T(x,y,t)=O(e^{-\lambda x})~~\hbox{for}~~x>0, \\
\label{eq:1.3a} |\nabla T|= O(e^{-\lambda|x|})~~
\hbox{for}~~\hbox{some}~~\lambda>0.
\end{eqnarray}
If such 
conditions are satisfied initially then they are valid for 
all subsequent times (see Section \ref{gub}). 

We are interested in a general situation,
when traveling waves solutions may not be relevant or not exist.

We introduce a natural quantity that measures the typical burning rate
\[
V(t)=\int\limits_{\Omega}T_t(x,y,t)\frac{dxdy}{H}.
\]
We call $V(t)$ the (instantaneous) ``bulk burning rate''
and its time average
\[
\langle V\rangle_t=\frac{1}{t}\int\limits_0^tV(s)ds,~~
\langle V\rangle_\infty=\liminf_{t\to\infty}\langle V\rangle_t
\]
simply ``bulk burning rate''. Because $T$ is non-dimensional,
$V$ has units of length per time, i.e. of velocity. 
Note that when $T(x,y,t)$ is a traveling wave front-like solution
 $T(x-ct,y,t)$, then
$V$ is the speed of the front, $V(t)=c$. But $V$ is defined for 
very general initial data and
more general equations and does not require assumptions about the
nature of 
the burning region. The word ``bulk'' refers to the fact that we take
a space (or space-time) average, capturing only bulk or large scale
effects.

Our first result shows that no matter 
what the advecting velocity is, it cannot slow down the bulk burning rate 
below a universal lower bound, of the same order of magnitude as the
laminar front speed.
Let us denote
\[ \alpha=-\inf_{0\le T\le 1}f'(T)>0, ~~\beta=-\sup_{0\le T\le 1}f''(T)>0. \]
(Note $\alpha = 1, \beta = 2$ for KPP).

\begin{thm}\label{thm2}
 There exists a constant $C>0$ such that for any
 advecting velocity $u(x,y,t)$ satisfying (\ref{eq:1.10}),
  (\ref{eq:1.11}) and (\ref{eq:1.11.1}),
  any solution of the passive-reactive diffusion
  equation (\ref{eq:1.1}) with  boundary conditions (\ref{eq:1.2})
  and (\ref{eq:1.3}) or (\ref{eq:1.2.1}), (\ref{eq:1.3}), the bulk burning 
rate $V(t)$ obeys the lower bound
  \begin{equation}
  \label{eq:1.5}
  V(t)\ge Cv_0\sqrt{\frac{\beta}{4\alpha}}(1-e^{-{\alpha v_0^2t}/{2\kappa}}).
  \end{equation}
\end{thm}
One of the applications of this theorem is 
in a homogenization regime, where the reaction is very weak
(see Appendix \ref{A1}).

As far as the general upper bound is concerned, 
it is easy to show (see Section \ref{gub}) for a very general class of velocities $u(x,y,t)$ that if the initial data $T_0(x,y)$ 
 satisfies (\ref{eq:1.3}) with $\lambda=v_0/ 2\kappa,$ then
\begin{equation}\label{eq:1.3.1}
\langle V\rangle_t\le\frac{L_0}{t}+\|u_1\|_{\infty}+v_0,
\end{equation}
with the constant length $L_0$ depending on the initial data $T_0$ only.
Here $\|u_1\|_\infty$ denotes, as usual, the supremum of $|u_1|$ over the whole
domain.
Therefore, the bulk burning rate may not exceed a linear bound in  the
amplitude of the advecting velocity.  For a
large class of flows  we prove  lower bounds on the bulk
burning rate that are 
linear in the magnitude of advection. For instance, a corollary
to Theorem 4 of  Section 4  
concerning mean zero shear flow of the form
\[ u(x,y)=(u(y),0), ~~\int_0^H u(y)dy=0 \]
can be stated simply as
\begin{thm}\label{thm1} There exists a constant $C>0$ that
depends only on the nonlinearity $f$ but not upon $u(y)$ and $T_0$
such that, for any solution  $T(x,y,t)$ of the 
  passive-reactive diffusion equation (\ref{eq:1.1}) with 
  boundary conditions (\ref{eq:1.2}), (\ref{eq:1.3}) or
  (\ref{eq:1.2.1}), (\ref{eq:1.3}) and any $\tau \geq \tau_0 = 
\hbox{max}[\frac{\kappa}{v_0^2}, \frac{H}{v_0}]$
the bulk burning rate obeys
\begin{equation}
  \label{eq:1.7}
  \langle V\rangle_\tau \ge C\left(1+\frac{\kappa^2}{v_0^2h_u^2}\right)^{-1}
\frac{\|u\|_{1}^2}{\|u\|_{\infty}}
\end{equation}
where $\displaystyle\|u\|_{1}=\int\limits_0^H|u(y)|\frac{dy}{H}$ and
$\displaystyle h_u=\frac{\|u\|_{1}}{\|u'\|_{\infty}}$ .
\end{thm}

Recall that the normalization  
of the reaction rate in (\ref{eq:1.1}) is chosen so that the laminar
traveling wave front speed is $v_0$ no matter what $\kappa$ is.
If we allow $\kappa$ to vary while keeping the coefficient $M$ in front of the
reaction term fixed, we find that the bound is still independent of
$\kappa,$ however the time $\tau_0$ in which it is reached behaves as
$\kappa^{-1/2}.$ 

Actually, we prove a far more general geometric 
estimate (Section \ref{shear}), 
from which (\ref{eq:1.7})
follows. This general estimate provides a non-trivial lower bound also
for the case when the
ratio of $L^1$ and $L^\infty$ norms of $u(y)$ becomes small. An
important feature of the bound (\ref{eq:1.7}) is the presence of the
ratio of the characteristic scale $h_u$ of variations of the advection
velocity and the reaction scale $l=\kappa/v_0$. The
estimate degenerates in the case $h_u\ll l$ which is expected from
physical considerations since additional wrinkling on the scales
smaller than the reaction scale should not accelerate the front.  

Furthermore, we consider time dependent shear flows, and obtain a 
similar lower bound on the bulk burning rate:
\[
\langle V\rangle\ge K(h_u,\tau_*,\frac{u}{|u|})\|u\|_{\infty}.
\]
Here $\tau_*$ is a typical time scale of the flow, defined similarly
(but not identically) to $h_u$. The prefactor $K$ becomes smaller when
either $h_u$ becomes smaller than the reaction scale $l=\kappa/v_0$,
or the time scale $\tau_*$ is faster than the reaction time
$\tau_c=\kappa/v_0^2$ or the time $\tau_H=H/v_0$ it takes the reaction to
traverse the cross-section. The precise formulation of our result
for the time dependent shear flows is given in Section \ref{timeshear}.

Finally, in Section \ref{perc} we consider a generalization of the time
 independent shear
flows. Namely, we consider ``percolating flows'',  
flows that have two or more (sufficiently regular) tubes of
streamlines connecting $x=-\infty$ and $x=+\infty$. 
These flows are not necessarily spatially periodic and can 
have completely arbitrary features outside the tubes of streamlines.  
We show that the bulk burning rate is still linear in the magnitude of the 
advecting velocity, no matter
what kind of behavior (closed streamlines, areas of still fluid, etc.) the flow
has outside the tubes. The proportionality coefficient depends on the 
geometry of 
the flow. Thus, we identify a broad class of the flows which increase the 
bulk burning rate linearly with the amplitude of the flow, the fastest possible
 rate of increase. 
We also  show that in general the dependence of the bulk burning rate
 $V(t)$ on the
magnitude of the advecting velocity may be sub-linear. An extreme
example is provided by shear 
flow perpendicular to the front (in periodic boundary conditions) 
where there is no significant enhancement: the bulk
burning rate remains uniformly
 bounded as the magnitude of the advecting velocity tends to infinity. 
Other examples are certain
cellular flows (flows
 with closed
streamlines): for every $\alpha\in(0,1)$ we construct cellular flows for
which $V(t)$ is bounded above by 
$C\|u\|_\infty^\alpha$ (Section \ref{examples}). 
A comparison of the results of this paper 
with extensive numerical studies will be presented in 
a companion paper with Fausto Cattaneo, Andrea Malagoli and Natalia Vladimirova
\cite{CCKMORV}.

For the rest of this paper, we assume that the reaction $f(T)$ satisfies 
(\ref{eq:1.4}), the initial data
satisfy (\ref{eq:1.3}), (\ref{eq:1.3a}),
and the advection velocity satisfies $\|u_1\|_\infty < \infty,$ 
$\| \nabla u \|_\infty <\infty.$
As we will see in the next section, such assumptions on $u$ ensure that 
the localization (\ref{eq:1.3}), (\ref{eq:1.3a})
of the initial data is preserved during the time evolution.
Throughout the paper we denote by $C$ various (not necessarily equal)
constants which may depend only on reaction $f(T)$.

\section{Preliminaries and an upper bound on the bulk burning rate}\label{gub}

Our considerations in this section follow the general ideas of \cite{Ber-Nir-2}. 
We show in this section that the boundary conditions (\ref{eq:1.3})
are conserved by evolution and establish the simple upper bound
(\ref{eq:1.3.1}).

\begin{lemma}\label{lemma0}
Assume that the initial data $T_0(x,y)$ satisfies the following bounds:
\begin{equation}
\label{eq1}
T(x,y)\le C_0e^{-\lambda x},~~
1-T(x,y)\le C_0e^{\lambda x},~~
|\nabla T|\le \frac{C_0}{H}e^{-\lambda |x|},~~
C_0>0,~\lambda>0.
\end{equation}
Let
$c_1\ge \|u_1\|_{\infty}+\kappa\lambda+\frac{v_0^2}{4\kappa\lambda}$,
  and $c_2\ge \|u_1\|_{\infty}+\kappa\lambda+\frac{v_0^2}{2\kappa\lambda}+
\frac{4\|\nabla u\|_\infty}
{\lambda}.$ Then
\begin{eqnarray}\label{eq:2.2.0}
&&T(t,x,y)\le C_0e^{-\lambda_0(x-c_1t)},~~~
1-T(x,y)\le C_0e^{\lambda_0(x+c_1t)}, \\
\label{eq:2.2.0.2}
&&|\nabla T|\le \frac{C_0}{H}e^{\mp\lambda(x\mp c_2t)}~~~for~~any~~t.
\end{eqnarray}
\end{lemma}
{\bf Proof.} The proof is an application of the maximum principle.
Note that $T$ satisfies an inequality
\[
  T_t+u\cdot\nabla T-\kappa\Delta T-\frac{v_0^2}{4\kappa}T\le 0.
\]
Introduce $\phi(x,t)=e^{-\lambda (x-c_1t)}$, then
\[
  \phi_t+u\cdot\nabla \phi-\kappa\Delta\phi-\frac{v_0^2}{4\kappa}\phi=
\lambda(c-u_1-\kappa\lambda-\frac{v_0^2}{4\kappa\lambda})\phi\ge 0
\]
by our assumptions on $c$. Applying the maximum principle to the function
$w=C_0\phi-T$ we obtain the first estimate in (\ref{eq:2.2.0}). To get 
the second estimate we note that $G=1-T$ satisfies the inequality
\[
G_t+u\cdot\nabla G-\kappa\Delta G\le 0.
\]
We let then $\psi(x,t)=e^{\lambda_0(x+c_1t)}$ and proceed as before
applying maximum principle to $w_1=C_0\psi-G$.

The decay of $|\nabla T|$ is obtained as follows. Let $P=|\nabla T|^2$,
then $P$ satisfies the equation
\[
  P_t+u\cdot\nabla P-\kappa\Delta P+2(|\nabla T_x|^2 +|\nabla T_y|^2)=
\frac{v_0^2}{2\kappa}f'(T)P-2(T_xu_x\cdot\nabla T+T_yu_y\cdot\nabla T).
\]
Therefore if we let $K=\frac{v_0^2}{2\kappa}+4\|\nabla u\|_{\infty}$,
then we get 
\[
   P_t+u\cdot\nabla P-\kappa\Delta P-KP\le 0
\]
and then (\ref{eq:2.2.0.2}) follows as before.$\Box$

Lemma \ref{lemma0} implies that the bulk burning rate cannot be larger
than
$c_1=v_0+\|u_1\|_{\infty}$ provided that the initial data decays fast
enough.
\begin{thm}\label{lemma0.1}
  Assume that $T_0(x,y)\le C_0e^{-\lambda x}$ and $1-T_0(x,y)\le
  C_0e^{\lambda x}$ with $\displaystyle \lambda=\frac{v_0}{2\kappa}$, then
  \begin{equation}
    \label{eq:2.2.6}
    \langle V\rangle_t \le \frac{4C_0\kappa}{v_0 t}+v_0+\|u_1\|_{\infty}.
  \end{equation}
\end{thm}
{\bf Proof.} We have
\[
\langle V\rangle (t)=\frac{1}{t}\int_0^tds\int\frac{dxdy}{H}T_s(s,x,y)=
\frac{1}{t}\int\frac{dxdy}{H}[T(x,y,t)-T_0(x,y)].
\]
Lemma \ref{lemma0} implies that $T(x,y,t)$ satisfies the bound
\[
T(x,y,t)\le C_0 e^{-\lambda(x-c_1t)},
\]
with $c_1=\|u_1\|_{\infty}+\kappa\lambda +\frac{v_0^2}{4\kappa\lambda}$.
Then we have
\begin{eqnarray}
  \label{eq:2.2.7}
\langle V\rangle (t)&\leq&\frac{1}{t}\int_0^H\frac{dy}{H}\int_{-\infty}^{0}
((1-T_0)-(1-T))\,dx+
 \frac{1}{t}\int_0^H\frac{dy}{H}\int_{0}^{c_1t}
[T(x,y,t)-T_0(x,y)]\nonumber
\\&+& \frac{1}{t}\int_0^H\frac{dy}{H}\int_{c_1t}^{\infty}
T\,dx\le\frac{2C_0}{t\lambda}+c_1\le 
\frac{4C_0\kappa}{t v_0}+\|u_1\|_\infty+v_0\nonumber
\end{eqnarray}
with our choice of $\lambda=\frac{v_0}{2\kappa}$.$\Box$

A simple corollary of Theorem \ref{lemma0.1} is that a shear flow in the direction
perpendicular to the front propagation does not enhance the bulk burning
rate.

\section{A universal lower bound for bulk burning rate}\label{glb}

In this section we prove Theorem \ref{thm2}. Let us integrate
(\ref{eq:1.1}) over the set $\Omega$ and obtain
\begin{equation}
  \label{eq:2.1}
  V(t)=\frac{v_0^2}{4\kappa}\int\limits_{\Omega}f(T)\frac{dxdy}{H}.
\end{equation}
Integration by parts is justified for any $t$ by Lemma \ref{lemma0}.
A straightforward computation that uses (\ref{eq:1.1}), the
boundary conditions (\ref{eq:1.2}) or (\ref{eq:1.2.1}), 
(\ref{eq:1.3}), (\ref{eq:1.3a}), and
the incompressibility of $u(x,y)$ shows that
\begin{equation}
  \label{eq:2.2}
  \frac{dV}{dt}\ge 
\frac{\beta v_0^2}{4}\int\limits_{\Omega}|\nabla T|^2\frac{dxdy}{H}-
\frac{\alpha v_0^2}{4\kappa}V.
\end{equation}
The two simple equations (\ref{eq:2.1}) and (\ref{eq:2.2}) are the basis of our technique for deriving
lower bounds on the burning velocity. 
The temperature drops from $1$ on the left to
$0$ on the right. The reaction $f(T)$ is large where $T$
takes values in some region strictly between $0$ and $1.$ If $T$ varies 
slowly, we get a good lower bound for $V(t)$
from (\ref{eq:2.1}). On the other hand, if $T$ varies fast, 
then the $|\nabla T|^2$ term in
(\ref{eq:2.2}) 
is large, which will give a lower bound on the bulk burning rate.
The equations 
(\ref{eq:2.1}) and (\ref{eq:2.2}) are thus complementary in this respect.
The following Lemma gives 
precise meaning to the statement that the reactive term in (\ref{eq:2.1})
and the gradient term in (\ref{eq:2.2}) cannot be simultaneously small: 
\begin{lemma}\label{lemma1}
  Let $f(T)$ be a function of concave KPP type (\ref{eq:1.4})
  and assume that the continuously differentiable
 function $T(x,y)$ satisfies the following
  assumptions:
\begin{itemize}
\item [(i)] $0\le T(x,y)\le 1$,
\item[(ii)] $\displaystyle\lim_{x\to -\infty}T(x,y)=1$, 
$\displaystyle\lim_{x\to
    +\infty}T(x,y)=0$ for every $y\in[0,1]$,
\end{itemize}
Then there exists a constant $C>0$ (independent of the function $T$) such that
\begin{equation}
  \label{eq:2.3}
  \int\limits_{\Omega}f(T)dxdy\int\limits_{\Omega}|\nabla T|^2dxdy\ge CH^2.
\end{equation}
\end{lemma}
{\bf Proof.} We may assume that $\int\limits_{\Omega}|\nabla T|^2dxdy<\infty$
and $\int\limits_{\Omega}f(T(x,y))dxdy<\infty$, otherwise (\ref{eq:2.3}) is
trivial. Let $y\in(0,1)$ be such that
\[
\int\limits_{-\infty}^\infty |\nabla T(x,y)|^2dx \leq
3\int\limits_{\Omega}|\nabla T|^2\frac{dxdy}{H}
\]
and
\[
\int\limits_{-\infty}^\infty f(T(x,y))dx\le 3\int\limits_{\Omega}f(T(x,y))\frac{dxdy}{H}.
\]
Then there exist $x_1$, $x_2$ such that $|T(x_1,y)-T(x_2,y)|\ge
1-\epsilon$ while $f(T(x,y))\ge C\epsilon$ for all $x\in(x_1,x_2)$
because of the boundary conditions (ii) and property (\ref{eq:1.4}) of
the non-linearity.  Then we have
\[
C\epsilon|x_1-x_2|\le 3\int\limits_{\Omega}f(T(x,y))\frac{dxdy}{H}
\]
and
\[
\frac{(1-\epsilon)^2}{|x_1-x_2|}\le3\int\limits_{\Omega}|\nabla T|^2\frac{dxdy}{H}.
\]
Multiplying these two equations we obtain
\[
\int\limits_{\Omega}f(T)dxdy\int\limits_{\Omega}|\nabla
T|^2dxdy\ge\frac{C\epsilon(1-\epsilon)^2H^2}{9},
\]
which proves Lemma \ref{lemma1}. $\Box$

Lemma \ref{lemma1}, (\ref{eq:2.1}), and (\ref{eq:2.2}) imply that 
\[
\frac{dV}{dt}+\frac{\alpha v_0^2}{4\kappa}V\ge C\frac{\beta v_0^4}{16\kappa V}.
\]
Therefore
\[
V^2(t)\ge\frac{C\beta v_0^2}{4\alpha}+e^{-\alpha v_0^2t/(2\kappa)}\left[V^2(0)-
\frac{C\beta v_0^2}{4\alpha}\right],
\]
from which Theorem \ref{thm2} follows. 

We remark that a simple variation of Lemma \ref{lemma1} allows to prove
Theorem \ref{thm2} for more general domains than a strip.  


\section{Bulk burning rate in shear flows}\label{shear}

\begin{figure}
\centerline{
\psfig{figure=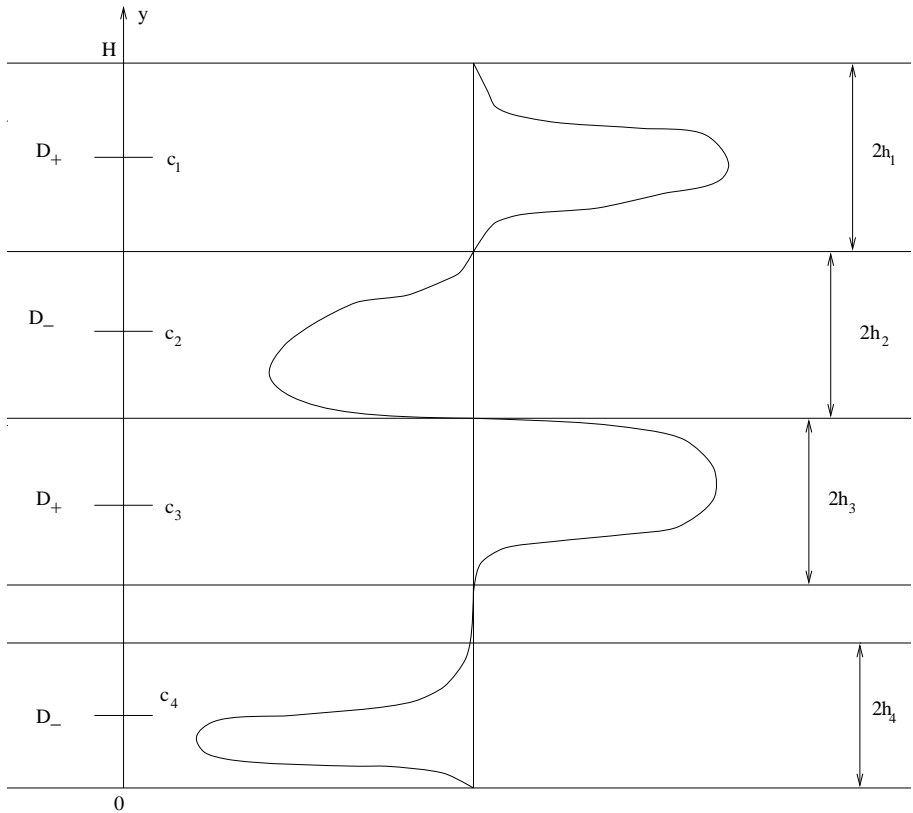,width=14cm } }

\caption{The structure of the shear flow.}
  \label{struct}
\end{figure}

Let us consider the passive-reactive diffusion
equation in a shear flow in two dimensions:
\begin{equation}
  \label{eq:3.1}
  T_t+u(y)T_x=\kappa\Delta T+\frac{v_0^2}{4\kappa}f(T).
\end{equation}
The boundary and initial
conditions are as in (\ref{eq:1.2}) or (\ref{eq:1.2.1}) and (\ref{eq:1.3}), (\ref{eq:1.3a}).
The flow $u(y)$ is continuously differentiable and has mean zero (a non-zero mean can
be taken into account by a simple change of variables):
\[
\int_0^H u(y)dy=0.
\]
We  prove now an estimate for the bulk burning
rate which is more general than the one presented in Theorem
\ref{thm1}. 
\begin{thm}\label{thm4}
  Let us consider an arbitrary partition of the interval $[0,H]$ into
  subintervals $I_j=[c_j-h_j,c_j+h_j]$ on which $u(y)$ does not change
  sign. Denote by $D_{-},$ $D_{+}$ the unions of intervals $I_j$ where $u(y)>0$
and $u(y)<0$ respectively (see figure \ref{struct}). 
Then there exists a constant $C>0$, independent of the
  partition, $u(y),$ and the initial data $T_0(x,y)$, so that the
  average burning rate $\langle V \rangle_\tau$
satisfies the
  following estimate:
\begin{equation}
  \label{eq:3.2}
  \langle V\rangle_\tau
 \ge C\left(
c_+
\sum_{I_j \subset D_+}
\left(1+\frac{l^2}{h_j^2}\right)^{-1}
\int\limits_{c_j-\frac{h_j}{2}}^{c_j+\frac{h_j}{2}}|u(y)|\frac{dy}{H} 
+c_-
\sum_{I_j \subset D_-}
\left(1+\frac{l^2}{h_j^2}\right)^{-1}
\int\limits_{c_j-\frac{h_j}{2}}^{c_j+\frac{h_j}{2}}|u(y)|\frac{dy}{H}\right) \nonumber
\end{equation}
for any $\tau \geq \tau_0=\hbox{max}\left[ \frac{\kappa}{v_0^2}, \frac{H}{v_0} \right].$
Here $l=\kappa/v_0.$ 
The constants $c_\pm$ are defined by 
\[ c_\pm = \left( \sum\limits_{I_j \subset D_\mp} \frac{h_j^3}{h_j^2 +l^2}
\right) \left( \sum\limits_{I_j} \frac{h_j^3}{h_j^2 +l^2}
\right)^{-1}. \]
\end{thm}

\it Remarks. \rm 1.  The exact choice of the unions of intervals 
$D_\pm$ is left to us. 
Given a shear flow, we should attempt to pick $D_\pm$ in a way to 
maximize the bound (\ref{eq:3.2}).
We provide a simple example after the proof  which 
illustrates how this bound works. \\
2. The proof of this bound becomes much easier, and also extends to more 
general types of reaction, if we assume $T_t(x,y,t) \geq 0.$ 
This is the case for traveling waves, the existence of which has been proven
for shear flows and various types of reaction in \cite{Ber-Nir-2}.
The minimal speed of traveling waves also provides
lower bounds for asymptotic speed of propagation of any front-like 
data \cite{Mal-Roq,jxin-2}. We plan to address the results which one can prove along
these lines in subsequent publications.
The disadvantage of an a priori assumption $T_t \geq 0,$ 
however, is that we cannot get an estimate on the time required to reach 
the lower bound from any initial data, and cannot extend the results to 
time-dependent flows, as we do in Section \ref{timeshear}. \\
 {\bf Proof.} 
  The plan of the proof is as follows. We know that 
 the bulk burning rate $V(t)$ satisfies the  inequalities
\begin{equation}\label{eq:3.2.1}
V(t)+\frac{4\kappa}{\alpha
  v_0^2}\frac{dV}{dt}\ge\frac{\beta\kappa}{\alpha}\int_{\Omega}|\nabla
T|^2dxdy.
\end{equation}
and 
\begin{equation} \label{eq:new1}
  V(t)=\frac{v_0^2}{4\kappa}\int_{\Omega}f(T)\frac{dxdy}{H}.
\end{equation}
We will be able to get a bound from below in terms of $u$ for the 
combination of the terms on the right-hand sides of these equations, which 
then provide a bound for $\langle V \rangle_\tau.$ 
The key will be to integrate in $x$ along the streamlines of the flow and 
then to average several times in $y$ and $t$ to bring $T_{yy}$ and $T_t$ terms
 to a form convenient for the estimate.

Consider one interval of $D_+$, $I_j=[c_j-h_j,c_j+h_j],$ and so that $u(y)>0$ for
$y\in(c_j-h_j,c_j+h_j)$. Integrating (\ref{eq:3.1}) over $x\in \Rm$,
we obtain 
\begin{equation}
  \label{eq:3.3}
\int\limits_{\Rm} T_t \,dx-\kappa \int\limits_{\Rm} T_{yy} \, dx -
\frac{v_0^2}{4\kappa}\int\limits_{\Rm}  f(T)\, dx=u(y). 
\end{equation}
Therefore
\begin{equation}
  \label{eq:3.3a}
\int\limits_{\Rm} T_t \,dx-\kappa\int\limits_{\Rm} T_{yy} \, dx \geq u(y). 
\end{equation}
(It might appear that we make the estimate rather crude by dropping  
the $f(T)$ term. There is however heuristic and numerical evidence
that in many situations
this term is insignificant in the regions where $u(y) > 0$,
and that the front is quite sharp in these regions. 
In contrast, in the regions where $u(y)<0$ 
the burning region is wider
and the $f(T)$ term will not be discarded there.)
Now we estimate both terms on the left hand side of (\ref{eq:3.3a}).
Let us begin with the second derivative term. To reduce the order of differentiation,
we employ the following averaging in $y:$
\begin{equation} \label{aver}
 \int\limits_0^{h_j/2} d \gamma \int\limits_{h_j/2-\gamma}^{h_j/2+\gamma}
d \delta \int\limits_{c_j-\delta}^{c_j+\delta} \cdot \,dy =
\int\limits_{c_j-h_j}^{c_j+h_j} G(h_j, y-c_j)\cdot \,dy, 
\end{equation}
where the kernel $G(h,\xi)$ can be computed explicitly as
\begin{eqnarray}\label{eq:defG}
G(h,\xi)=\left\{\begin{matrix}
{\frac{1}{2}(h-|\xi|)^2-\left(\frac{h}{2}-|\xi|\right)^2, & |\xi|<h/2\cr
\frac{1}{2}(h-|\xi|)^2, &h/2 \leq |\xi|<h\cr}
\end{matrix}\right.
\end{eqnarray}
Observe that the function $G(h,\xi)$ has the following properties
\begin{eqnarray}
  \label{eq:3.11}
  &&G(h,\xi)\ge 0~~\hbox{for all}~~\xi\in[-h,h] \nonumber\\
&&G(h,\xi)\le \frac{h^2}{4}~~\hbox{for}~~
\xi\in [-h,h]  \\
&& G(h, \xi) \geq \frac{h^2}{8} ~~\hbox{for}~~
\xi\in[-\frac{h}{2},\frac{h}{2}]. \nonumber
\end{eqnarray}
Let us apply the averaging procedure (\ref{aver}) to the equation
(\ref{eq:3.3a}). We have
\begin{lemma} \label{Tyy}
\begin{equation} \label{eqn1}
\left| \int\limits_{\Rm} dx \int\limits_{c_j-h_j}^{c_j+h_j} G(h_j, y-c_j)
T_{yy}(x,y)\,dy \right| \leq C\left(  h_j^2 
\int\limits_{\Rm} dx \int\limits_{c_j-h_j}^{c_j+h_j} |\nabla T|^2 \,dy +
 \int\limits_{\Rm} dx \int\limits_{c_j-h_j}^{c_j+h_j} f(T) \,dy \right).
\end{equation}
\end{lemma}
{\bf Proof.} We are going to split $x \in \Rm$ into two sets. In one set, the 
$L^2_y$ norm of the gradient will be large and we will use it to estimate
the second derivative term. In the other set, the variation of the 
temperature will be small, and we will use the reactive term to bound the
second derivative term. More precisely, let $\rho$ be a number such 
that $\displaystyle\sqrt{2\rho}h_j=1/3$
and define the set ${\cal D}_{j\rho}\subset \Rm$ by
\[
{\cal D}_{j \rho}=\left\{x \in \Rm:\int\limits_{c_j-h_j}^{c_j+h_j}|\nabla T|^2(x,y)dy\ge\rho
  h_j\right\}
\]
so that
\begin{equation}
  \label{eq:3.6.0}
\int\limits_{c_j-h_j}^{c_j+h_j}dy|T_y(x,y)|\le
\sqrt{\frac{2}{\rho}}\int\limits_{c_j-h_j}^{c_j+h_j}|\nabla T(x,y)|^2dy=
6h_j\int\limits_{c_j-h_j}^{c_j+h_j}|\nabla T(x,y)|^2dy
\end{equation}
for $x\in{\cal D}_{j\rho}$.
Notice that for such $x,$ according to (\ref{aver}),
\begin{eqnarray} \label{eqn2}
 && \int\limits_{c_j-h_j}^{c_j+h_j} G(h_j, y-c_j)
T_{yy}(x,y)\,dy =  \int\limits_0^{h_j/2} d \gamma
 \int\limits_{h_j/2-\gamma}^{h_j/2+\gamma} d \delta 
(T_y(x,c_j+\delta)-T_y(x,c_j-\delta)) \\ && \leq 
3h_j^2 \int\limits_{c_j-h_j}^{c_j+h_j}|\nabla T(x,y)|^2dy. \nonumber
\end{eqnarray}
For $x$ outside ${\cal D}_{j\rho},$ we use the representation
\begin{equation} \label{eqn4}
   \int\limits_{c_j-h_j}^{c_j+h_j} G(h_j, y-c_j)
T_{yy}(x,y)\,dy=
\int\limits_{h_j/2<|y-c_j|\leq h_j}T(x,y)dy 
-\int\limits_{|y-c_j|\leq h_j/2}T(x,y)dy. 
\end{equation}
We need the following crucial observation:
\begin{lemma}\label{lemma2}
Assume that 
\[ |T(x,y_1)-T(x,y_2)| \leq \frac{1}{3}, \]
 then we have
\[
|T(x,y_1)-T(x,y_2)|\le C(f(T(x,y_1)+f(T(x,y_2))
\]
for any $y_1,y_2\in(c_j-h_j,c_j+h_j)$.
\end{lemma}
{\bf Proof.} Let us denote $T_1=T(x,y_1)$,
$T_2=T(x,y_2)$, then we have using (\ref{eq:1.4})
\[
f(T_1)+f(T_2)\ge \inf_{T\in(|T_1-T_2|,1-|T_1-T_2|)}f(T)\ge C|T_1-T_2|.
\]
This proves Lemma \ref{lemma2}. $\Box$ \\
 Notice 
 that if $x\notin{\cal D}_{j \rho}$ we have
\[ 
\int\limits_{c_j-h_j}^{c_j+h_j}|T_y(x,y)|dy\le\sqrt{2\rho}h=\frac{1}{3}
\]
and therefore $|T(x,y_1)-T(x,y_2)|< 1/3$ for any
$y_1,y_2\in(c_j-h_j,c_j+h_j)$.
Applying Lemma \ref{lemma2} to (\ref{eqn4}), we get 
\[ \left| \int\limits_{\Rm} dx \int\limits_{c_j-h_j}^{c_j+h_j} G(h_j, y-c_j)
T_{yy}(x,y)\,dy \right| \leq C  
\int\limits_{\Rm} dx \int\limits_{c_j-h_j}^{c_j+h_j} f(T) \, dy. \]
This completes the proof of Lemma \ref{Tyy}. $\Box$ \\
Averaging equation (\ref{eq:3.3a}) and applying Lemma \ref{Tyy} we
obtain
\begin{eqnarray}
  \label{eq:3.8}
  &&\int\limits_{c_j-h_j}^{c_j+h_j}dyu(y)G(h_j,y-c_j)\le
\int\limits_{\Rm} dx\int\limits_{c_j-h_j}^{c_j+h_j}dy G(h_j,y-c_j)T_t(x,y)\\&&
+C \left(h_j^2\kappa \int\limits_{\Rm} dx
\int\limits_{c_j-h_j}^{c_j+h_j}dy |\nabla T(x,y)|^2
+\kappa
\int\limits_{\Rm} dx \int\limits_{c_j-h_j}^{c_j+h_j}dy f(T(x,y))\right).\nonumber
\end{eqnarray}

Next we consider the intervals where velocity $u(y)\le 0$. These are
analyzed similarly, except that now we do not discard the last term on
the right side in (\ref{eq:3.3}). The estimate analogous to
(\ref{eq:3.8}) is
\begin{eqnarray}
  \label{eq:3.9}
 &&\int\limits_{c_j-h_j}^{c_j+h_j}dy|u(y)|G(h_j,y-c_j)\le
\frac{v_0^2}{4\kappa} 
\int\limits_{\Rm} dx\int\limits_{c_j-h_j}^{c_j+h_j}dy G(h_j,y-c_j)f(T(x,y)) \\&&
-\int\limits_{\Rm} dx\int\limits_{c_j-h_j}^{c_j+h_j}dy G(h_j,y-c_j)T_t(x,y)
\nonumber \\&&
+C \left(h_j^2\kappa \int\limits_{\Rm} dx
\int\limits_{c_j-h_j}^{c_j+h_j}dy |\nabla T(x,y)|^2
+\kappa
\int\limits_{\Rm} dx \int\limits_{c_j-h_j}^{c_j+h_j}dy f(T(x,y))\right).\nonumber
\end{eqnarray}

Thus we succeeded in replacing the second order derivative term with 
expressions directly linked to the burning rate. Now it remains to estimate
the time derivative term.
Summing (\ref{eq:3.8}) and (\ref{eq:3.9}) over $I_j \subset D_+$ and 
$I_j \subset D_-$ respectively, and using properties (\ref{eq:3.11})
of the kernel $G,$ we get
\begin{eqnarray}
  \label{eq:3.13}
  &&\kappa\int\limits_{D_+}dy\int\limits_{\Rm} dx|\nabla T(x,y,t)|^2+
\frac{v_0^2}{4\kappa}\int\limits_{D_+}dy\int\limits_{\Rm} dxf(T(t,x,y))\\&&+
\sum_{I_j\subset D_+}
\int\limits_{I_j}dy\int\limits_{\Rm}dx\frac{G(h_j,y-c_j)}{h_j^2+l^2}T_t(x,y,t)
\ge C\sum_{I_j\subset D_+}
\left(1+\frac{l^2}{h_j^2}\right)^{-1}
\int\limits_{c_j-\frac{h_j}{2}}^{c_j+\frac{h_j}{2}}|u(y)|dy\nonumber
\end{eqnarray}
and
\begin{eqnarray}
  \label{eq:3.14}
  &&\kappa\int\limits_{D_-}dy\int\limits_{\Rm} dx|\nabla T(x,y,t)|^2+
\frac{v_0^2}{4\kappa}\int\limits_{D_-}dy\int\limits_{\Rm} dxf(T(t,x,y))\\&&-
\sum_{I_j\subset D_-}
\int\limits_{I_j}dy\int\limits_{\Rm} dx\frac{G(h_j,y-c_j)}{h_j^2+l^2}T_t(x,y,t)
\ge C\sum_{I_j\subset D_-}
\left(1+\frac{l^2}{h_j^2}\right)^{-1}
\int\limits_{c_j-\frac{h_j}{2}}^{c_j+\frac{h_j}{2}}|u(y)|dy\nonumber
\end{eqnarray}
(recall $l=\kappa/v_0$).
Let us choose the weights $m_+$ and $m_-$ according to
\begin{equation}\label{eq:3.141}
m_\pm = \sum_{I_j\subset D_\mp}
\int_{I_j}dy\frac{G(h_j,y-c_j)}{h_j^2+l^2}
\end{equation}
Set also $M= \hbox{max}(m_+, m_-)$.
Notice that by the properties of $G,$ we have 
\begin{equation}\label{csms}
 \frac{1}{16} c_\pm \leq \frac{m_\pm}{M} \leq \frac{1}{4}c_\pm 
\end{equation}
for constants $c_\pm$ in the formulation of the Theorem.
Let us define measures
\begin{equation}\label{eqn6}
d \nu_\pm =\sum\limits_{I_j \subset D_\pm}\frac{m_\pm \chi_{I_j}(y) G(h_j, y-c_j)}{M H (h_j^2+l^2)}dy. 
\end{equation}
Here we denote, as usual, by $\chi_S$ the characteristic function
of the set $S.$ 
Multiplying (\ref{eq:3.13}) and (\ref{eq:3.14}) by $m_+$ and
$m_-$ respectively, and adding them together we obtain
\begin{eqnarray} \label{tog}
\kappa \int\limits_\Omega |\nabla T|^2 \frac{dxdy}{H} +
\frac{v_0^2}{4\kappa} \int\limits_\Omega f(T) \frac{dxdy}{H} +
\int\limits_{\Rm} dx \int\limits_0^H d\nu_+(y) T_t -
\int\limits_{\Rm} dx \int\limits_0^H d\nu_-(y) T_t \\
\geq C \left( \frac{m_+}{M}\sum_{I_j\subset D_+}
\left(1+\frac{l^2}{h_j^2}\right)^{-1}
\int\limits_{c_j-\frac{h_j}{2}}^{c_j+\frac{h_j}{2}}|u(y)|dy
+ \frac{m_-}{M}\sum_{I_j\subset D_-}
\left(1+\frac{l^2}{h_j^2}\right)^{-1}
\int\limits_{c_j-\frac{h_j}{2}}^{c_j+\frac{h_j}{2}}|u(y)|dy
\right). \nonumber
\end{eqnarray}
We have the following
\begin{lemma}\label{Tt}
For any $\tau_1,$ $\tau_2,$ 
\begin{eqnarray*}
\int\limits_{\tau_1}^{\tau_2}dt \left(\,
 \int\limits_{\Rm}dx \int\limits_0^H d\nu_+(y) T_t -
\int\limits_{\Rm} dx \int\limits_0^H d\nu_-(y) T_t \right) \\
\leq 
C \sum_{i=1}^2
\left( \frac{H\kappa}{v_0} \int\limits_\Omega |\nabla T(x,y,\tau_i)|^2
\frac{dxdy}{H} 
+ \left(1+\frac{Hv_0}{\kappa} \right) \int\limits_\Omega f(T(x,y,\tau_i)) 
\frac{dxdy}{H} \right).
\end{eqnarray*}
\end{lemma}
{\bf Proof.} By definition (\ref{eqn6}) of the measures $\nu_\pm,$ their total weights
are equal:
\begin{equation}\label{eqne1}
\int\limits_0^H d\nu_+(y) = \int\limits_0^H d \nu_-(y) <1. 
\end{equation}
 It is easy to construct a measure preserving bijective map 
$\Phi(y): D_+ \rightarrow D_-,$ so that $\nu_+(S) = \nu_-(\Phi(S)).$ 
Then we can write 
\begin{eqnarray}\label{eqn5}
 \int\limits_{\tau_1}^{\tau_2}dt \left(
 \int\limits_{\Rm} dx \int\limits_0^H d\nu_+(y) T_t -
\int\limits_{\Rm} dx \int\limits_0^H d\nu_-(y) T_t \right) \\ 
= \int\limits_{\Rm} dx \int\limits_0^H d\nu_+(y)(T(x, \Phi(y), \tau_1)-
T(x,y,\tau_1)) +  \int\limits_{\Rm} dx \int\limits_0^H d\nu_+(y)(T(x,y, \tau_2)-
T(x,\Phi(y),\tau_2)). \nonumber
\end{eqnarray}
Consider the first term on the left hand side of (\ref{eqn5}).
We split all $x \in \Rm$ into two sets, saying $x \in S$ if there
exists $y$ such that 
\[ |T(x,y,\tau_1)- T(x,\Phi(y),\tau_1)| > \frac{1}{3}. \]
Using the same argument  we applied in the proof of Lemma \ref{lemma1}, one can show 
that
\[ \int\limits_0^H |\nabla T(x,y)|^2 \,dy \int\limits_0^H f(T(x,y))\,dy \geq C \left[ \hbox{sup}_{y_1,y_2 \in [0,H]}
|T(x,y_1)-T(x,y_2)| \right]^3, \]
where $C$ is some universal constant, depending only on $f.$ Therefore,
for every $x \in S,$ 
 \[
\left( \int\limits_0^H |\nabla T(x,y)|^2 \,dy \int\limits_0^H f(T(x,y)) \right)^{\frac{1}{2}} \geq C, 
\]
and hence, using the fact that the total weight of $\nu_+$ does not exceed $1$ by (\ref{eqne1}), we have
\begin{equation} \label{S1}
 \int\limits_0^H d\nu_+(y)|(T(x, \Phi(y), \tau_1)-
T(x,y,\tau_1))| \leq C \left( \frac{\kappa}{v_0} \int\limits_0^H |\nabla T(x,y)|^2 \, dy + \frac{v_0}{\kappa}
\int\limits_0^H f(T(x,y))\, dy \right). 
\end{equation}
For $x \notin S,$ we have 
\begin{equation}\label{S2}
 \int\limits_0^H d \nu_+(y)  | T(x,y,\tau_1) - T(x, \Phi(y), \tau_1)| \leq C  \int\limits_0^H d \nu_+(y)[f(T(x,y)) +
f(T(x,\Phi(y)))]  \leq C \int\limits_0^H f(T(x,y))\frac{dy}{H} 
\end{equation}
by Lemma \ref{lemma2} and (\ref{eqn6}).
Equations (\ref{S1}) and (\ref{eqn6}) together imply Lemma \ref{Tt}. $\Box$ \\
Theorem \ref{thm4} now follows from Lemma \ref{Tt}, relations
(\ref{eq:3.2.1}) and (\ref{eq:new1}), and inequality (\ref{tog}).
Given time interval $[0,\tau],$ apply the following averaging to the
both sides of (\ref{tog}):
\begin{equation}\label{avt}
  \frac{1}{\tau^3} \int\limits_0^{\frac{\tau}{4}} d\gamma \int\limits_
{\frac{\tau}{4}-\gamma}^{\frac{\tau}{4}+\gamma} d \delta \int\limits_{\frac{\tau}{2}-\delta}
^{\frac{\tau}{2}+\delta} dt = \frac{1}{\tau^3} \int\limits_0^\tau G(\frac{\tau}{2}, t-\frac{\tau}{2})\, dt. 
\end{equation} 
A direct computation using (\ref{eq:3.2.1}) and (\ref{eq:new1}) shows that the left 
hand side of (\ref{tog}) after averaging  does not exceed
\begin{equation}\label{fin}
  \frac{C}{\tau} \int\limits_0^\tau V(t)\,dt \left( 1 + \frac{H/v_0 + \kappa /v_0^2}{\tau} +
\frac{H \kappa}{v_0^3 \tau^2} \right). 
\end{equation}
On the other hand, the right hand side of (\ref{tog}) is independent of time, and averaging (\ref{avt}) results
in multiplication by a constant independent of $\tau.$ 
Hence from (\ref{fin}) we see that if 
\[ \tau \geq \tau_0=\hbox{max}[\frac{\kappa}{v_0^2}, \frac{H}{v_0}], \] 
then  
\[ \langle V \rangle_\tau \geq  C \left( \frac{m_+}{M}\sum_{I_j\subset D_+}
\left(1+\frac{l^2}{h_j^2}\right)^{-1}
\int\limits_{c_j-\frac{h_j}{2}}^{c_j+\frac{h_j}{2}}|u(y)|dy
+ \frac{m_-}{M}\sum_{I_j\subset D_-}
\left(1+\frac{l^2}{h_j^2}\right)^{-1}
\int\limits_{c_j-\frac{h_j}{2}}^{c_j+\frac{h_j}{2}}|u(y)|dy
\right). \]
This proves Theorem \ref{thm4} (recall that $m_\pm$ and $c_\pm$ are related by (\ref{csms})).
$\Box$ \\

\it Remark. \rm  At the expense of making the proof slightly more technical, 
 the total width of the strip $H$ in the formula for $\tau_0$
can be replaced by an often smaller value $\tilde{H},$ which is introduced as follows.
Consider the function $g(y) = \nu_+[0,y] -\nu_{-}[0,y]$ on $[0,H].$ 
Then $\tilde{H}$ is defined as a maximal distance between two neighboring roots
of $g.$ 
It is straightforward to generalize the proof of Lemma \ref{Tt} to yield this result,
by taking a specific measure preserving function $\Phi$ which maps $D_+$ to $D_-$ only within 
the intervals between the neighboring roots of $g.$ 
The characteristic time $\frac{\tilde{H}}{v_0}$ has a clear intuitive physical meaning:
this is the time needed to burn across the scale on which the shear flow $u$ wrinkles
 the front. \\

\it Example. \rm  Let $u(y)=u_0 \sin \frac{2\pi n y}{H}.$ We can take intervals $I_j$ as half-periods
of $u$ where it does not change sign. Factors $c_\pm$ are equal in this case.
Then Theorem \ref{thm4} implies
that for any $\tau \geq \tau_0 = \hbox{max}[ \frac{\kappa}{v_0^2}, \frac{H}{v_0}]$ we have
\[ \langle V \rangle_\tau \geq C\left(1+ \frac{n^2l^2}{H^2}\right)^{-1}u_0. \]
According to the above Remark, it is easy to see that in this example, 
$H$ in definition of $\tau_0$ can be replaced with $\tilde{H}=H/n.$
The map $\Phi$ in this example can be taken to map half-periods where
$u$ is positive on the neighboring half-periods where $u$ is negative. \\

\noindent Theorem \ref{thm1} is a direct corollary of Theorem \ref{thm4}. \\
{\bf Proof.} Let us denote by $|S|$  the Lebesgue measure of set $S.$ 
Define the sets
\[
{\cal F_\pm}=\left\{y\in[0,H]:~\pm u(y) \ge  \frac{1}{4}
\int\limits_0^H|u(y)|\frac{dy}{H}=\frac{1}{4}\|u\|_{1}\right\}.
\]
Notice that since $u$ is mean zero, $\|u_\pm\|_1 = \frac{1}{2}\|u\|_1$
(here $u_\pm$ are the positive and the negative parts of $u$). 
Therefore
\[ \frac{1}{4} \|u\|_{1}(H-|{\cal F}_\pm|) + \|u\|_{\infty}
|{\cal F}_\pm| \geq \frac{1}{2}\|u\|_{1}H, \]
and so 
\begin{equation}\label{Fpm}
 |{\cal F}_\pm|  \geq \frac{\|u\|_{1}H}{4\|u\|_{\infty}}. 
\end{equation}
Let $\displaystyle
h_u=\frac{\|u\|_{L^1}}{\|u'\|_{L^\infty}}$, then for any
$y\in{\cal F_\pm}$ and any $y'\in(y-h_u/8,y+h_u/8)$ we have $\displaystyle
|u(y')|\ge \frac{1}{8}\|u\|_{L^1}$. 
It is easy to construct unions
$\displaystyle D_\pm = \cup_j I_{j}^{\pm}$ of non-overlapping intervals
$I^{\pm}_{j} = (y_{j}^{\pm}-h_u/8, y_{j}^{\pm}+h_u/8)$ with $y_{j}^{\pm}
 \in \cal F_\pm,$ such that $\displaystyle 
|D_\pm|  \geq \frac{1}{2} |{\cal F}_\pm|.$ Then we have:
\begin{equation}\label{Dpm}
\int\limits_{D_\pm}|u(y)|\frac{dy}{H}\ge \frac{1}{16}\frac{\|u\|_{1}}{H}
|{\cal F}_\pm|.
\end{equation}
Combining (\ref{Fpm}), (\ref{Dpm}), and (\ref{eq:3.2}) we get
\[ \langle V \rangle_\tau  \geq C\left(1+\frac{l^2}{ h^2_u} \right)^{-1} 
 \frac{\|u\|_1^2}{\|u\|_\infty}. \,\,\,\,\,\,\,\,\,\,\,\,\,\Box \]

\section{Time dependent shear flows}\label{timeshear}

One may expect linear growth of the bulk burning rate  in the amplitude
of the flow, but the temporal characteristic scale $\tau_*$ of the 
variations of the flow will play a role similar to that of the scale $h_u$ in
Theorem \ref{thm1}. That is, too rapid oscillations will diminish the
enhancement of the bulk burning rate. 
Consider two systems of intervals $I_j$ in $[0,H],$ $D_+$ and $D_-.$
At this point we do not make any assumptions regarding the
behavior of $u(y,t)$ on
$D_\pm.$ Let us introduce the notation
\begin{eqnarray}
  \label{eq:50.1}
 J(t,u)& =& c_+ \sum_{I_j\in D_+}
\left(1+\frac{l^2}{h_j^2}\right)^{-1}
\int\limits_{c_j-\frac{h_j}{2}}^{c_j+\frac{h_j}{2}}u(t,y)\frac{dy}{H} \\&& -
c_-\sum_{I_j\in D_-}
\left(1+\frac{l^2}{h_j^2}\right)^{-1}
\int\limits_{c_j-\frac{h_j}{2}}^{c_j+\frac{h_j}{2}}u(t,y)\frac{dy}{H},\nonumber
\end{eqnarray}
where $c_\pm$ are defined as in Theorem \ref{thm4}. 
Given a starting time $t_0$ and length of the time interval $\tau,$ 
we define 
\begin{equation}\label{eqn10}
 J(t_0, \tau, u) = \frac{1}{\tau^3} \int\limits_{t_0}^{t_0+\tau}
G(\frac{\tau}{2}, t-t_0-\frac{\tau}{2}) J(t,u)\,dt.
\end{equation}
We also denote $\langle V \rangle_{t_0, \tau}$ the average of
$V(t)$ over  an interval of time of duration $\tau$ starting at time $t_0:$ 
\[ \langle V \rangle_{t_0, \tau} = \frac{1}{\tau}
\int\limits_{t_0}^{t_0+\tau} V(t) \, dt. \]
We have
\begin{thm}\label{thm4.1}
For any choice of the intervals $I_j \subset [0,H]$ in $D_\pm$ and
any $\tau,$ $t_0$ 
\begin{equation}
  \label{eq:50.1.1}
  \langle V\rangle_{t_0,\tau} \ge C
\left(1+\left(
\frac{\tau_0}{\tau}\right)^2\right)^{-1} J(t_0, \tau, u),
\end{equation}
where $\tau_0= \hbox{max}[\frac{\kappa}{v_0^2}, \frac{H}{v_0}].$ 
\end{thm}
{\bf Proof.} The proof is a direct corollary of the proof of Theorem 
\ref{thm4}.
We remark that the choice of the averaging in time (\ref{avt}) in the last
stage of the proof is not the only one possible. Given a particular time
dependent flow at some initial moment $t_0,$
one can try to adjust the averaging procedure to get a better lower bound.
 $\Box$

\it Remark. \rm Similarly to the remark after the proof of Theorem 
\ref{thm4}, we can replace $H$ in the definition of $\tau_0$ by a smaller
value $\tilde{H}$ (defined in that remark). \\

To clarify the meaning of Theorem \ref{thm4.1} we make several observations
and consider two examples. The general way one can apply this 
theorem is as follows. Given a moment of time $t_0,$ we try to choose
$\tau$ and $D_\pm$ so as to maximize the lower bound (\ref{eq:50.1.1}).
There is a certain tradeoff involved in choosing $\tau.$ If we take $\tau$
to be small, it is likely that we can find $D_\pm$ so that velocity 
$u(y,t)$ does not change sign there during time interval $[t_0, t_0+\tau],$ 
staying positive on $D_+$ and negative on $D_-.$ Then there is no cancellation
in equation (\ref{eqn10}) defining $J(t_0, \tau, u).$ 
However, the factor $(1+\frac{\tau_0^2}{\tau^2})^{-1}$ may become very small
if $\tau \ll \tau_0.$ If we take $\tau$ large, some cancellation is 
likely to occur in (\ref{eqn10}), making the bound weaker,
 unless the shear flow  varies
on time scales still larger than $\tau.$ 
In the flows which oscillate in time on the scale smaller than $\tau_0$, 
we will not be able to avoid either cancellation in  (\ref{eqn10}) or 
small factor $(\tau/\tau_0)^2$ in the bound  (\ref{eq:50.1.1}), and will 
end up with weaker lower bound than we would have gotten if the flow varied
slower in time. Notice that in any case the bound grows linearly 
with the amplitude of the flow, and there are factors reflecting the 
moderating effect of fast oscillations both in space and in time.

If we want to know the average of the bulk burning rate
over a long period of time, much larger than the typical time scale of
the flow, we can use Theorem \ref{thm4.1} by splitting this long time period
into appropriately chosen 
smaller ones and getting lower bounds on the averages over these
smaller time intervals. 
Combined, they will also give us an estimate on the long time average. \\
 
\it Example 1. \rm 
Consider a flow
\[ u(y,t) = u_0 \sin 2\pi \omega t \sin \frac{2\pi n y}{H}. \]
To get an estimate on long-time average of the burning rate for 
such flow, consider $t_0=0.$ Set $\tau = \frac{1}{2\omega},$ and 
take 
\[ D_+ = \bigcup\limits_{j=1}^{n}I_{2j-1}, \,\,\,D_- =
\bigcup\limits_{j=1}^n I_{2j}, \]
where $I_j= (\frac{(j-1)H}{2n}, \frac{jH}{2n}).$
Then we get 
\begin{eqnarray*}
J(0,\tau,u)&=& u_0(1+4 (\tau_0\omega)^2
)^{-1} \left(1+\frac{n^2 l^2}{H^2}
\right)^{-1} \\ && \times
\int\limits_0^\tau dt \frac{G(\tau/2, t-\tau/2)}{\tau^3}
\sin 2\pi \omega t \sum\limits_{j=1}^{2n}
\int\limits_{\frac{(j-3/4)H}{2n}}^{\frac{(j-1/4)H}{2n}} 
|\sin \frac{2\pi n y}{H}|\,dy \\ & \geq &
C(1+4 (\tau_0 \omega)^2)^{-1} \left(1+\frac{n^2 l^2}{H^2}
\right)^{-1}u_0,
\end{eqnarray*}
with constant $C$ dependent only on reaction $f$.
A similar estimate is valid for $t_0=\frac{1}{2\omega},$
we only need to switch $D_\pm.$ Therefore, we get that 
for any $t_0$ and any $\tau_1 \geq \frac{1}{\omega},$
\begin{equation}\label{ex1}
 \langle V \rangle_{t_0,\tau_1} \geq C(1+4 (\tau_0 \omega)^2)^{-1} 
\left(1+\frac{n^2 l^2}{H^2} \right)^{-1}u_0.
\end{equation}
It is not difficult to obtain estimates for averages
over times smaller than $\frac{1}{\omega},$ but these would
generally (and naturally) depend on the choice of starting 
time $t_0.$  \\
\it Example 2. \rm 
Consider 
\[ u(y,t) = u_0 \sin \frac{2\pi n (y-ct)}{H}. \]
This is a flow which shifts in $y$ direction. We assume for 
simplicity that the boundary conditions for $T$ are periodic in $y.$
Given any $t_0,$ pick $\tau= \frac{H}{8cn}.$ Time $\tau$ is chosen so that 
during this time, the regions where $u$ is positive and negative do 
not shift completely; there are regions where velocity stays 
positive or negative during $[0, \tau].$ Take 
\[ D_+ = \bigcup\limits_{j=1}^n I_{2j-1}, \,\,\,\,
D_-= \bigcup\limits_{j=1}^n I_{2j-1} \]
with $I_j = (\frac{(4j-3)H}{8n}, \frac{(4j-1)H}{8nH}).$ 
A direct computation of $J(t_0, \tau, u)$ shows the following 
bound:
\begin{equation}\label{ex2}
 \langle V \rangle_{t_0, \tau} \geq C\left(1+
\left(\frac{8cn\tau_0}{H}\right)^2\right)^{-1}
\left(1+\frac{n^2 l^2}{H^2} \right)^{-1}u_0, 
\end{equation}
where $C$ may depend only on reaction function $f.$ 
We note that in this example, it is to easy to show that  
(\ref{ex2}) extends to any averaging time $\tau,$ independently 
of the starting time $t_0:$
\[  \langle V \rangle_{t_0, \tau} \geq C\left(1+\frac{\tau_0^2}{\tau^2}\right)^{-1}
\left(1+\frac{n^2 l^2}{H^2} \right)^{-1}u_0. \] 

\section{Percolating flows}\label{perc}

We now consider a more general class of flows, which we call ``percolating''. 
By this we mean  that there exist at least two tubes of streamlines
of the advecting velocity $u(x,y)$, one of which connects $x=-\infty$
and $x=+\infty$, and the other one goes from $x=+\infty$ to
$x=-\infty$. More precisely, let us assume that there exist
regions $D_j^+$ and $D_j^-$, $j=1,\dots N$ such that each of them is
bounded by the streamlines of $u(x,y)$, and the projection of each
streamline of $u(x,y)$, contained in either $D_j^+$ or $D_j^-$, onto
the $x$-axis covers the whole real line (these projections need not be
one-to-one, however). As before, we denote $D_\pm$ the union of all
$D_j^\pm$ respectively.

Our considerations in this section will follow closely the ideas
of the shear flow case.
However, there are two natural geometries in the problem. The 
Laplace operator is
best described in Euclidean coordinates, while for the advection term the 
geometry of streamlines imposed by the flow is most natural.
In the case of the shear flows these geometries coincide, but generally
they are at odds. Due to this fact, additional technical difficulties
arise when we consider percolating flows.

\begin{figure}
  \centerline{
  \psfig{figure=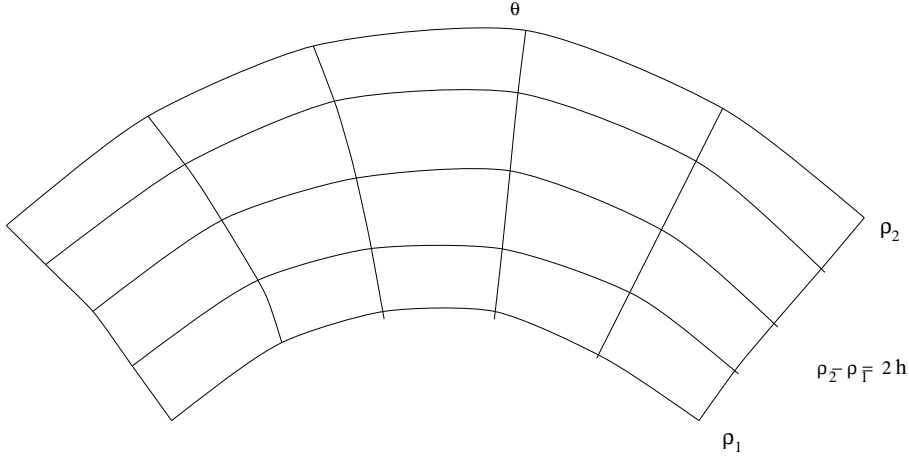,width=12cm} 
     }
  \caption{Curvilinear coordinates $(\rho,\theta)$.}
  \label{fig-coord}
\end{figure}

We assume that the streamlines in $D_j^\pm$
are sufficiently regular, so that
inside each $D_j^\pm$ there exists a one-to-one $C^2$ change of
coordinates $(x,y) \rightarrow (\rho, \theta),$ such that $\rho$ is
constant on the streamlines, while $\theta$ is an orthogonal
coordinate for $\rho$ (with a slight abuse of notation we shall use the
same notation $(\rho,\theta)$ in all $D_j^\pm,$ although these coordinates are
not defined globally).
Moreover,
$u\cdot\nabla\theta>0$ in $D_j^+$, while $u\cdot\nabla\theta<0$ in
each $D_j^-$. On $D_j^\pm,$ $\rho$ varies in $[c_j^\pm-h_j^\pm,
c_j^\pm+h_j^\pm],$ while $\theta$ varies in $(-\infty, \infty).$ See
figure \ref{fig-coord} for a sketch of coordinates $(\rho, \theta).$ The
square of the length element inside each set $D_j^{\pm}$ is given by
\[
dx^2+dy^2=E_1^2(\rho,\theta)d\rho^2+E_2^2(\rho,\theta)d\theta^2.
\]
 
We assume that the functions $E_{1,2}$ satisfy the following
conditions. They are bounded from above and below:
\begin{equation}
  \label{eq:5.0}
  C^{-1}\le E_{1,2}(\rho,\theta)\le C
\end{equation}
uniformly on all $D_j^\pm.$ Moreover, the function 
\begin{equation}\label{eq:5.00}
\omega(\rho,\theta)=\frac{E_2(\rho,\theta)}{E_1(\rho,\theta)}
\end{equation}
satisfies the following bounds:
\begin{equation}
  \label{eq:5.1}
  C^{-1}\le|\omega(\rho,\theta)|\le C,~~
\left|\frac{\partial\omega}{\partial\rho}(\rho,\theta)\right|
\le\frac{C}{h_j^{\pm}} \,\,\,\hbox{on}\,\,\,D_j^\pm \,\,\,\hbox{respectively},
\end{equation}
with $2h_j^\pm$ being the absolute value of the difference of the
values of $\rho$ on the two components of the boundary $\partial
D_j^{\pm}$ (recall that $D_j^{\pm}$ are bounded by two streamlines of
$u(x,y)$). 

Conditions (\ref{eq:5.0}) and (\ref{eq:5.1}) are satisfied, for
instance, in the following examples:
\begin{itemize}
\item[(1)]{ A flow $u(x,y)=U\tilde u(x,y)$ with $U$ being a scalar, and the
    flow $\tilde u(x,y)$ satisfying on $D_j^\pm$
    \begin{equation}\label{ucon}
      C^{-1}\le|\tilde u(x,y)|\le C,~~|\nabla\times\tilde u(x,y)|\le 
\frac{C}{d},
    \end{equation}
where $d$ is the maximum length of a level set of $\theta$ inside $D_j^\pm$.}
\item[(2)] {More generally, it is enough to ask that in each $D_j^\pm$ 
there exists 
a function $\psi$ constant on the streamlines of $u$ such that 
\begin{equation}\label{psi}
C^{-1} \leq |\nabla \psi| \leq C, \,\,\,|\Delta \psi|\leq \frac{C}{d}.
\end{equation}}
\end{itemize}
We remark that (1) is a particular case of (2) where $\psi$ is taken 
to be a stream function of the flow $\tilde u$. 

\begin{figure}
  \centerline{
  \psfig{figure=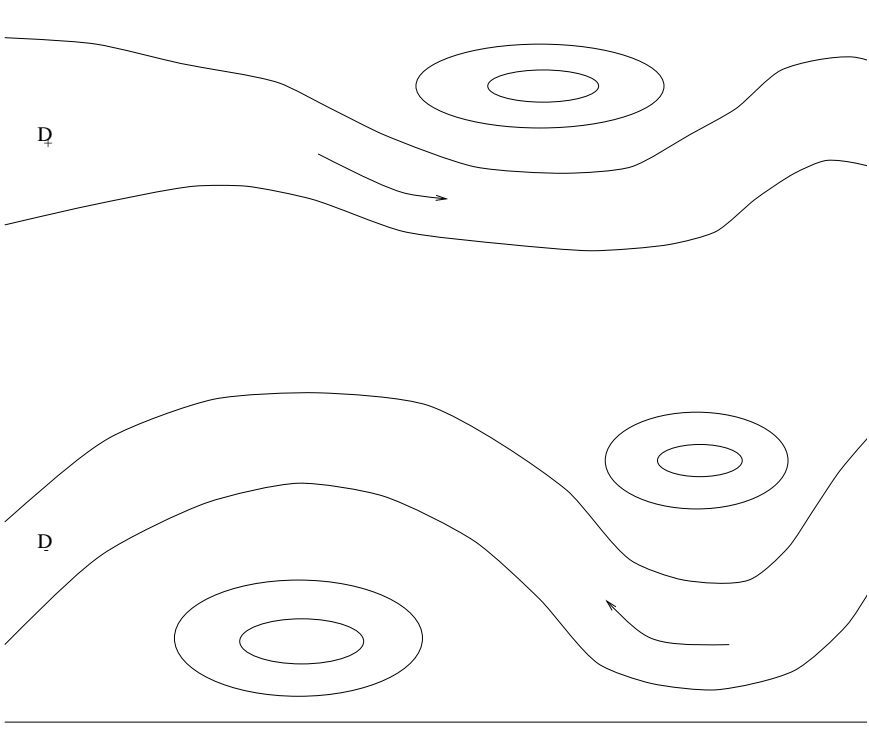,width=14cm} 
     }
  \caption{Streamlines of $u(x,y)$.}
  \label{fig1}
\end{figure}

We do not make any assumptions on the behavior of the
streamlines of $u(x,y)$ outside the regions $D_+$ and $D_-$. In
particular, there may be pockets of still fluid, streamlines may be
closed, etc. (see Figure \ref{fig1}).
 
Another assumption concerns the relative measure of the sets $D_\pm$,
which we assume to oscillate not too wildly. More precisely, let
$R_{ab}$ denote the rectangle $R_{ab}=[a,b]_x\times[0,H]_y$, and let
$D_\pm^{ab}=D_\pm\cap R_{ab}$. We define the measures
\[
\mu_\pm[a,b]=\int\limits_{D_\pm^{ab}}d\rho d\theta
\sum_{j}E_1(\rho,\theta)E_2(\rho,\theta)\frac{G(h_j,\rho-c_j)}{H(h_j^2+l^2)}
\]
with $l=\kappa/v_0$ and the function $G(h,\rho)$ defined by
(\ref{eq:defG}). Then we assume that there exists a partition of the
real axis
\begin{equation}
  \label{eq:parti}
  \dots<x_{-n}<\dots<x_{-1}<x_0<x_1<\dots x_n<\dots, \,\,\,x_{i+1}-x_i \leq L
\,\,\,\hbox{for all}\,\,\,i \in Z,
\end{equation}
and a number $m_0$ so that the ratio
\begin{equation}
  \label{eq:measratio}
  m_0=\frac{\mu_+[x_i,x_{i+1}]}{\mu_-[x_i,x_{i+1}]}~~~
\hbox{is independent of $i\in Z$}.
\end{equation}
This assumption is not the weakest necessary assumption and we make it
for clarity of exposition. The simplest example where this assumption is satisfied
is periodic percolating flows for which (\ref{ucon}) or (\ref{psi})
is satisfied. 

 Then we have the following Theorem.
\begin{thm}\label{thm5}
  Let each of the sets $D_j^{\pm}$ be of the form
  $D_j^{\pm}=\left\{\rho\in[c_j-h_j,c_j+h_j]\right\}$. Then under the
    assumptions made above, we have
\begin{eqnarray}
  \label{eq:5.2}
  \langle V\rangle_\tau&\ge& C\left(\frac{1}{1+m_0}\sum_{D_j^+}
\left(1+\frac{l^2}{h_j^2}\right)^{-1}
\int\limits_{c_j-\frac{h_j}{2}}^{c_j+\frac{h_j}{2}}|u(\rho, \theta)|
E_1(\rho, \theta)
\frac{d\rho}{H} \right. \\
&+& \frac{1}{1+m_0^{-1}}\sum_{D_j^-}\left.
\left(1+\frac{l^2}{h_j^2}\right)^{-1}
\int\limits_{c_j-\frac{h_j}{2}}^{c_j+\frac{h_j}{2}}|u(\rho, \theta)|
E_1(\rho, \theta)\frac{d\rho}{H}\right)
\nonumber
\end{eqnarray}
for every 
\[ \tau \geq \tau_0 = \hbox{max}\left[\frac{\kappa}{v_0^2}, 
\frac{H+L}{v_0}\right]. \]
Here $l=\kappa/v_0$, $L$ and $m_0$ are as in (\ref{eq:parti})
and (\ref{eq:measratio}), 
and the constant $C$ in (\ref{eq:5.2}) depends
on the function $f(T)$ and the constants appearing in (\ref{eq:5.0})
and (\ref{eq:5.1}).
\end{thm}
\it Remark. \rm Notice that the integrals on 
the right-hand side are independent of $\theta$
and give fluxes of the fluid through the middles of the tubes of 
streamlines.

\noindent {\bf Proof.} The proof of this theorem follows the steps of the proof
of Theorem \ref{thm4} for the shear flow. We will again utilize the
differential inequality (\ref{eq:2.2}), and the expression (\ref{eq:2.1}) for
the bulk burning rate, as well as multiple averaging over regions bounded by the
streamlines of the advecting velocity $u(x,y)$. 

Let us consider one region
$D_j^+=\left\{(\rho,\theta):~\rho\in[c_j-h_j,c_j+h_j]\right\}$.
 Let us also denote the tube of streamlines
$D_\delta=\left\{(\rho,\theta):~\rho\in[c_j-\delta,c_j+\delta],
  \theta\in(-\infty,\infty)\right\}$, $\delta<h_j$ and integrate (\ref{eq:1.1})
over the set $D_\delta\subset D_j^+$:
\begin{equation}
  \label{eq:5.5}
  \int\limits_{D_\delta}T_tdxdy-\int\limits_{-\delta}^\delta
d\rho u(\rho,\theta)E_1(\rho,\theta)+\kappa\int\limits_{D_\delta}\Delta Tdxdy=
\frac{v_0^2}{4\kappa}\int\limits_{D_{\delta}}f(T)dxdy.
\end{equation}
We used here the relation
\begin{equation}\label{eq:5.6}
\int\limits_{D_\delta}u\cdot\nabla Tdxdy=-\int\limits_{-\delta}^\delta
u(\rho,\theta)E_1(\rho,\theta,)d\rho
\end{equation}   
which follows from the fact that $D_\delta$ is the tube of
streamlines, and the boundary conditions (\ref{eq:1.3}). Moreover, the
quantity on the right side of (\ref{eq:5.6}) is independent of
$\theta$ because $u(x,y)$ is incompressible (\ref{eq:1.10}). 

We first estimate the term that involves the Laplacian in
(\ref{eq:5.5}). The following analog of Lemma \ref{Tyy} holds:
\begin{lemma} \label{Laplacian}
\[ \left|  \int\limits_0^{h_j/2}d\gamma
\int\limits_{h_j/2-\gamma}^{h_j/2+\gamma}d\delta \int\limits_{D_\delta}
\Delta T(x,y)dxdy \right| \leq C \int\limits_{D_j^+}
[ h_j^2 |\nabla T|^2 +f(T)]\,dxdy. \]
\end{lemma}
{\bf Proof}. Notice that 
\begin{eqnarray}
  \label{eq:5.7}
  \int\limits_{D_\delta}
\Delta T(x,y)dxdy&=&\int_{-\infty}^\infty d\theta\left[
\frac{E_2}{E_1}\frac{\partial T}{\partial\rho}(\delta,\theta)-
\frac{E_2}{E_1}\frac{\partial T}{\partial \rho}(-\delta,\theta)\right]\\
&=&\int_{-\infty}^\infty d\theta\left[\omega(\delta,\theta)
\frac{\partial T}{\partial\rho}(\delta,\theta)-\omega(-\delta,\theta)
\frac{\partial T}{\partial\rho}(-\delta,\theta)\right]\nonumber
\end{eqnarray}
with $\omega(\theta,\rho)$ defined in (\ref{eq:5.00}). Next, following
the general procedure in the proof of Theorem \ref{thm4} we fix
$\theta\in \Rm$ and average (\ref{eq:5.7}) over
$\delta\in[h_j/2-\gamma,h_j/2+\gamma]$ with $\gamma\in[0,h_j/2]$, and then
also average in $\gamma$. Then (\ref{eq:5.7}) becomes
\begin{equation}
  \label{eq:5.8}
  \int\limits_0^{h_j/2}d\gamma\int\limits_{h_j/2-\gamma}^{h_j/2+\gamma}d\delta
\left[T_{\rho}(\delta,\theta)\omega(\delta,\theta)-
T_\rho(-\delta,\theta)\omega(-\delta,\theta)\right].
\end{equation}
We show how to estimate the first term in (\ref{eq:5.8}), with the second
term treated in the same way. We integrate it by parts to get
\begin{eqnarray}
&& \int\limits_0^{h_j/2}d\gamma
\int\limits_{h_j/2-\gamma}^{h_j/2+\gamma}d\delta
T_{\rho}(\delta,\theta)\omega(\delta,\theta)=
\int\limits_{0}^{h_j/2}d\gamma\left[
T(\frac {h_j}2+\gamma,\theta)\omega(\frac {h_j}2+\gamma,\theta)\right.
\nonumber\\
&&\left.-
T(\frac{h_j}2-\gamma,\theta)\omega(\frac {h_j}2-\gamma,\theta)\right]-
\int\limits_0^{h_j/2}d\gamma\int\limits_{h_j/2-\gamma}^{h_j/2+\gamma}d\delta
T(\delta,\theta)\omega_\rho(\delta,\theta)\nonumber\\&&\label{eq:5.9}=
\int\limits_{0}^{h_j/2}
d\gamma\left[\left(T(\frac {h_j}2+\gamma,\theta)-1\right)
\omega(\frac {h_j}2+\gamma,\theta)-
\left(T(\frac {h_j}2-\gamma,\theta)-1\right)\omega(\frac {h_j}2-\gamma,\theta)
\right]\\&&-
\int\limits_0^{h_j/2}d\gamma\int\limits_{h_j/2-\gamma}^{h_j/2+\gamma}d\delta
\left(T(\delta,\theta)-1\right)\omega_\rho(\delta,\theta).\nonumber
\end{eqnarray}
Consider the set of $\theta$ such that
\[
\int\limits_{0}^{h_j}|T_\rho(\rho,\theta)|^2d\rho\ge \frac 1{4h_j}.
\]
We have for such $\theta$
\begin{equation}
  \label{eq:5.10}
 \int\limits_0^{h_j/2}d\gamma\int\limits_{h_j/2-\gamma}^{h_j/2+\gamma}d\delta
|T_\rho(\delta,\theta)\omega(\delta,\theta) |\le C
\int\limits_0^{h_j/2}d\gamma\sqrt{2\gamma}\left(\int\limits_0^{h_j}d\delta 
T_\rho^2(\delta,\theta)\right)^{1/2}\le C
h_j^2\int\limits_0^{h_j}T_\rho^2d\rho.
\end{equation}
Next we look at $\theta$ such that 
\[
\int\limits_{0}^{h_j}|T_\rho|^2d\rho\le \frac{1}{4h_j}.
\]
In this case for any $\rho_1,\rho_2\in[0,h_j]$ we have
\[
|T(\rho_1,\theta)-T(\rho_2,\theta)|\le\sqrt{h_j}
\left(\int\limits_{0}^{h_j}T_\rho^2d\rho\right)^{1/2}\le \frac 12.
\]
Therefore, either $T(\rho,\theta)\ge 1/4$ for all $\rho$, or
$|T(\rho,\theta)-1|\ge 1/4$ for all $\rho$. Then we have, $1-T \le
Cf(T)$, or $T\le Cf(T)$, respectively. We use one of these bounds and
the corresponding part of (\ref{eq:5.9}) to get for such $\theta$
\begin{eqnarray}
&& \left|\int\limits_0^{h_j/2}d\gamma
\int\limits_{h_j/2-\gamma}^{h_j/2+\gamma}d\delta
T_{\rho}(\delta,\theta)\omega(\delta,\theta)\right|\le
C\int\limits_{0}^{h_j}f(T(\rho,\theta))d\rho+
\frac{C}{h_j}
\int\limits_0^{h_j/2}d\gamma\int\limits_{h_j/2-\gamma}^{h_j/2+\gamma}
d\delta f(T(\delta,\theta))
\nonumber\\&&
\label{eq:5.11}\le C\int\limits_{0}^{h_j}f(T(\rho,\theta))d\rho.
\end{eqnarray}
Now we put together the estimates (\ref{eq:5.10}) and (\ref{eq:5.11})
to obtain for all $\theta\in \Rm$:
\begin{eqnarray}
  \label{eq:5.12}
  &&\left|\int\limits_0^{h_j/2}d\gamma\int\limits_{h_j/2-\gamma}^{h_j/2+\gamma}d\delta
\int\limits_{-\infty}^\infty d\theta
T_{\rho}(\delta,\theta)\omega(\delta,\theta)\right|\\&&
\le
C\left[
h_j^2\int\limits_{-\infty}^\infty d\theta\int\limits_{0}^{h_j}d\rho T_\rho^2+
\int\limits_{-\infty}^\infty d\theta\int\limits_{0}^{h_j}
d\rho f(T(\rho,\theta))\right]\nonumber
\\&&
\le  
C\int\limits_{D_j^+}\left[h_j^2|\nabla T|^2+f(T(x,y))\right]dxdy. \,\,\,\,\,\,\,\Box
\nonumber
\end{eqnarray}

Similarly to the shear case, 
Lemma \ref{Laplacian} and (\ref{eq:5.5}) imply the inequality
\begin{eqnarray}
&&-\int\limits_{D_j^+}
\frac{d\rho d\theta}{H} E_1(\rho,\theta)E_2(\rho,\theta)
\frac{G(h_j,\rho_j-c_j)}{(h_j^2+l^2)}
T_t(\rho,\theta)+
\int\limits_{c_j-h_j}^{c_j+h_j}\frac{d\rho}{H} 
\frac{G(h_j,\rho-c_j)}{h_j^2+l^2}
|u(\rho,\theta)|E_1(\rho,\theta)\nonumber\\
&&\le C\left[\kappa\int\limits_{D_j^+}|\nabla T|^2\,\frac{dxdy}{H}+
\frac{v_0^2}{4\kappa}
\int\limits_{D_j^+}
f(T(x,y))\,\frac{dxdy}{H} \right] \label{eq:5.15}
\end{eqnarray}
where $G(h, \xi)$ is defined as before by (\ref{eq:defG}).
An estimate similar to (\ref{eq:5.15}) holds in the regions $D_j^-$,
where the flow is going backwards, except that the time
derivative term in (\ref{eq:5.15}) enters now with the opposite sign:
\begin{eqnarray}
  \label{eq:5.15.1}
&&\int\limits_{D_j^-}
\frac{d\rho d\theta}{H} E_1(\rho,\theta)E_2(\rho,\theta)
\frac{G(h_j,\rho_j-c_j)}{h_j^2+l^2}
T_t(\rho,\theta)+
\int\limits_{c_j-h_j}^{c_j+h_j}\frac{d\rho}{H} 
\frac{G(h_j,\rho-c_j)}{h_j^2+l^2}
|u(\rho,\theta)|E_1(\rho,\theta)\\
&&\le C\left[\kappa\int\limits_{D_j^-}|\nabla T|^2\,
\frac{dxdy}{H}+\frac{v_0^2}{4\kappa}
\int\limits_{D_j^-}
f(T(x,y))\,\frac{dxdy}{H} \right].\nonumber  
\end{eqnarray}
Let us choose the weights 
\[
m_+=\frac{1}{1+m_0},~~~ m_-=\frac{1}{1+m_0^{-1}}, 
\]
(recall $m_0$ is defined by (\ref{eq:measratio})) so that
\begin{eqnarray}
  \label{eq:5.16}
  m_+ \mu_+ [x_i,x_{i+1}]=
 m_- \mu_- [x_i, x_{i+1}]
\end{eqnarray}
for any two points $x_i,$ $x_{i+1}$ of the partition (\ref{eq:parti})
of the $x$-axis, similarly to what we did in the proof of Theorem
\ref{thm4} (see (\ref{eq:3.141})).  In order to finish the proof of
Theorem \ref{thm5} we multiply equations (\ref{eq:5.15}) and
(\ref{eq:5.15.1}) by $m_+$ and $m_-$, respectively, and add them.
It remains now to estimate the time derivative term, and the following 
general Lemma provides us with the analog of Lemma \ref{Tt} in the shear 
case.
\begin{lemma}\label{lemma3}
  Let $\Omega_0$ be a rectangle $\Omega_0=L\times H$, and let
  $\Omega_{1,2}\subset\Omega_0$ be two open subsets of $\Omega_0$.
  Consider two continuous non-negative functions
  $\phi_{1,2}:\Omega_{1,2}\to \Rm$ such that $0 \le \phi_{1,2}(x,y)\le C$ and
\begin{eqnarray}
  \label{eq:5.18}
  \int\limits_{\Omega_1}dxdy\phi_1(x,y)=\int\limits_{\Omega_2}dxdy\phi_2(x,y).
\end{eqnarray}
Let $T:\Omega_0\to \Rm$ be a continuously differentiable function, $0\le
T\le 1$, then for any $\varepsilon>0$ we have
\begin{eqnarray}
  \label{eq:5.19}
 &&\left|\int\limits_{\Omega_1}dxdy\phi_1(x,y)T(x,y)
-\int\limits_{\Omega_2}dxdy\phi_2(x,y)
T(x,y)\right|\\&&\le C
\left[(L+H)\left[\varepsilon \int\limits_{\Omega_0}dxdy|\nabla T|^2+
\frac 1{\varepsilon}
\int\limits_{\Omega_0}f(T(x,y))dxdy\right]+\int\limits_{\Omega_0}
f(T(x,y))dxdy\right]. \nonumber
\end{eqnarray}
\end{lemma}
We postpone the proof of Lemma \ref{lemma3} till the end of this section. 
Using Lemma \ref{lemma3} in each rectangle $R_{x_i,x_{i+1}}$ with
$\varepsilon=\frac{\kappa}{v_0}$, $\Omega_1=D^+\cap
R_{x_i,x_{i+1}}$, $\Omega_2=D^-\cap
R_{x_i,x_{i+1}},$ and functions $\phi_{1,2}$ given by
\begin{eqnarray}
&&\phi_1=m_+\sum_{D_j^+}
\frac{G(h_j,\rho-c_j)}{H(h_j^2+l^2)}E_1(\rho,\theta)E_2(\rho,\theta)
\chi_{D_j^+}(\rho,\theta)\nonumber\\
&&\phi_2=m_-\sum_{D_j^-}
\frac{G(h_j,\rho-c_j)}{H(h_j^2+l^2)}E_1(\rho,\theta)E_2(\rho,\theta)
\chi_{D_j^-}(\rho,\theta),\nonumber
\end{eqnarray}
we arrive at the analog of Lemma \ref{Tt}:
\begin{eqnarray}\label{Ttp}
\int\limits_{\tau_1}^{\tau_2} dt\left( \int\limits_{D_j^+}
m_+\sum_{D_j^+}
\frac{G(h_j,\rho-c_j)}{H(h_j^2+l^2)}E_1(\rho,\theta)E_2(\rho,\theta)
\chi_{D_j^+}(\rho, \theta)T_t(\rho, \theta) d \rho d\theta \right. \\
-\left. \int\limits_{D_j^-}
m_-\sum_{D_j^+}
\frac{G(h_j,\rho-c_j)}{H(h_j^2+l^2)}E_1(\rho,\theta)E_2(\rho,\theta)
\chi_{D_j^-}(\rho, \theta)T_t(\rho, \theta) d \rho d\theta \right)\nonumber \\
\leq C  \sum_{i=1}^2
\left( \frac{(H+L)\kappa}{v_0} \int\limits_\Omega |\nabla T(x,y,\tau_i)|^2
\frac{dxdy}{H} 
+ \left(1+\frac{(H+L)v_0}{\kappa} \right) \int\limits_\Omega f(T(x,y,\tau_i)) 
\frac{dxdy}{H} \right). \nonumber
\end{eqnarray}
The rest of the proof of Theorem \ref{thm5} is completely analogous to the 
proof of Theorem \ref{thm4}. We average in time according to (\ref{avt})
and use (\ref{eq:2.1}) and (\ref{eq:2.2}) to conclude the proof. $\Box$

\begin{figure}
  \centerline{
  \psfig{figure=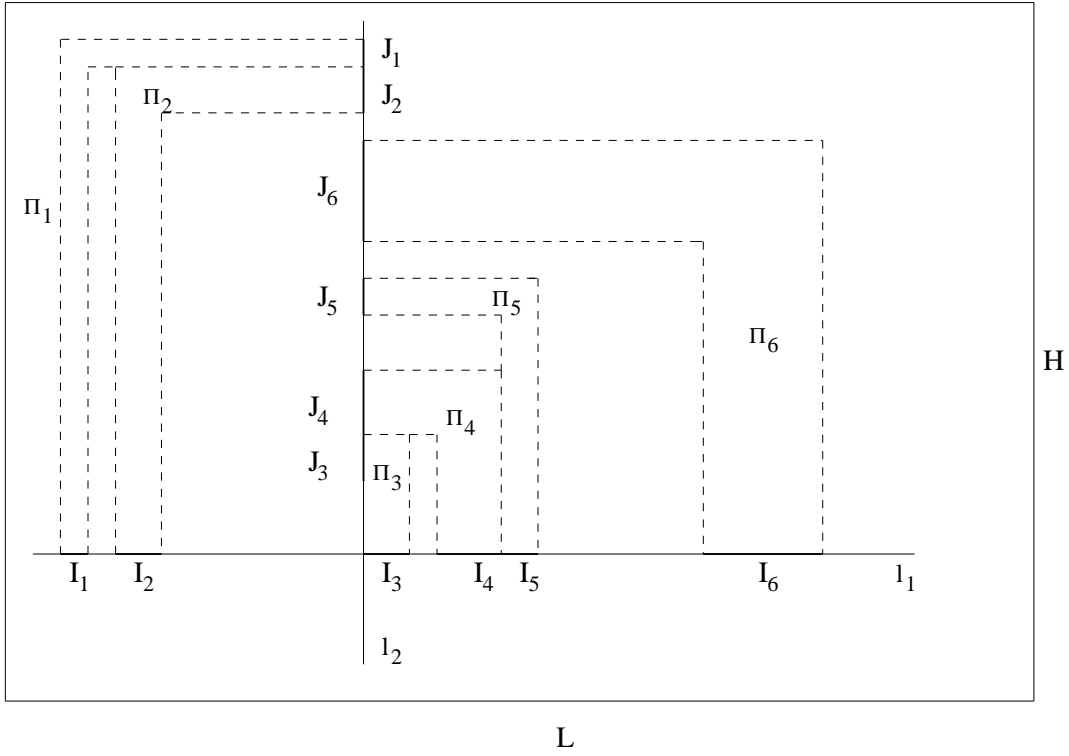,height=10cm} 
     }
  \caption{Staircase for $I_k$ and $J_k$.}
  \label{stair}
\end{figure}

We now give the proof of Lemma \ref{lemma3}. \\
{\bf Proof.} We define the measures
$\nu_{1,2}$ by
\[
d\nu_{1,2}(x,y)=\phi_{1,2}(x,y)\chi_{\Omega_{1,2}}(x,y)dxdy.
\]
Let ${\cal A}\subset\Omega_1$
be the set of points where $T(x,y)>5/8$, and let the open set ${\cal
  B}\subset\Omega_2$ be such that $\nu_1({\cal A})=\nu_2({\cal B})$.
Then we have
\begin{eqnarray}
  \label{eq:lemma.1}
  &&\int\limits_{\Omega_1}d\nu_1(x,y) T(x,y)-
\int\limits_{\Omega_2}d \nu_2(x,y) T(x,y)\le 
C\int\limits_{\Omega_1\backslash {\cal A}}dxdyf(T(x,y))\\&&+
\int\limits_{\cal A}d \nu_1(x,y)T(x,y)-
\int\limits_{\cal B}d \nu_2(x,y) T(x,y).\nonumber
\end{eqnarray}
Let us decompose further ${\cal B}={\cal B}'\cup{\cal B}''$, where
${\cal B}'=\left\{(x,y)\in{\cal B}:~T(x,y)>3/4\right\}$. We also
consider an open set ${\cal A}'\subset{\cal A}$ such that $\nu_1({\cal
  A}')=\nu_2({\cal B}')$, and write ${\cal A}={\cal A}'\cup{\cal
  A}''$. Then we obtain
\begin{eqnarray}
\int\limits_{\cal A}d nu_1 T-
\int\limits_{\cal B}d\nu_2 T=
\int\limits_{{\cal A}'}d \nu_1 T-
\int\limits_{{\cal B}'}d \nu_2 T+
\int\limits_{{\cal A}''}d \nu_1 T-
\int\limits_{{\cal B}''}d \nu_2 T\label{eq:lemma.2}
\end{eqnarray}
and, moreover, Lemma \ref{lemma2} implies that
\begin{eqnarray}
  \label{eq:lemma.3}
  \int\limits_{{\cal A}'}d\nu_1 T(x,y)-
\int\limits_{{\cal B}'}d\nu_2 T(x,y)\le C
\int\limits_{{\cal A}'}dxdyf(T(x,y))+
\int\limits_{{\cal B}'}dxdyf(T(x,y)).
\end{eqnarray}
Therefore, we are done if $\nu_1({\cal A}'')=\nu_2({\cal B}'')=0$.
Assume now that this is not the case. Then we may find a horizontal
line $l_1:~y=y_0$ and a vertical line $l_2:~x=x_0$ such that
\[
|l_1\cap{\cal A}''|\ge\frac{C}{L}\nu_1({\cal A}''),~~
|l_2\cap{\cal B}''|\ge\frac{C}{H}\nu_2({\cal B}''),
\]
where $|S|$ denotes the one-dimensional Lebesgue measure. Moreover, we
may choose subsets $Q_1\subset l_1\cap{\cal A}''$ and $Q_2\subset
l_2\cap{\cal B}''$ so that $Q_1=\cup_{k=1}^NI_k$ and
$Q_2=\cup_{k=1}^MJ_k$ are finite unions of intervals, and
$\displaystyle |Q_1|=|Q_2|\ge \frac{C\nu_1({\cal A}'')}{L+H}$. 
We may assume (possibly after subdividing into smaller intervals) that
$N=M$, and $|I_k|=|J_k|$ for all $k$. Let us connect each pair of
intervals $I_k$ and $J_k$ by perpendicular lines ``staircase''
as depicted on 
Figure \ref{stair}. Notice that for every point $(x,y) \in I_k,$
$T(x,y) >3/4$ while for every point $(x',y') \in J_k,$
$T(x',y') < 5/8.$
An argument directly analogous to that in the proof of Lemma
\ref{lemma1} shows that the following estimate holds:
\[
\int\limits_{\Pi_k}f(T)dxdy\int\limits_{\Pi_k}|\nabla T|^2dxdy\ge C|I_k|^2.
\]
Therefore we have for any $\varepsilon>0$
\begin{eqnarray}\label{eq:lemma.4}
&&\frac{1}{\varepsilon}\int\limits_{\Omega_0}f(T)dxdy+\varepsilon
\int\limits_{\Omega_0}|\nabla T|^2dxdy\ge \frac{C}{L+H}\nu_1({\cal A}'')\\&&\ge
\frac{C}{L+H}\left[~~\int\limits_{{\cal A}''}dxdy\phi_1(x,y)T(x,y)-
\int\limits_{{\cal B}''}dxdy\phi_2(x,y)T(x,y)\right].\nonumber
\end{eqnarray}
Equations (\ref{eq:lemma.1}-\ref{eq:lemma.4}) show that
\begin{eqnarray}
  \label{eq:lemma.5}
 &&\int\limits_{\Omega_1}dxdy\phi_1(x,y)T(x,y)-\int\limits_{\Omega_2}dxdy\phi_2(x,y)
T(x,y)\\&&\le C
\left[(L+H)\left[\varepsilon \int\limits_{\Omega_0}dxdy|\nabla T|^2+
\frac{1}{\varepsilon}
\int\limits_{\Omega_0}f(T(x,y))dxdy\right]+\int\limits_{\Omega_0}
f(T(x,y))dxdy\right]. \nonumber
\end{eqnarray}
The same proof shows that this bound holds for 
$\int_{\Omega_2}dxdy\phi_2(x,y)T(x,y)-\int_{\Omega_1}dxdy\phi_1(x,y)T(x,y)$.
This finishes the proof of Lemma \ref{lemma3}. $\Box$

\section{Examples of sub-linear growth of the bulk burning rate}
\label{examples} 

We give in this section examples of flows for which bulk burning rate grows
sub-linearly in the amplitude of the advecting velocity. We do not try
to identify the most general class of such flows but 
consider rather one simple family of flows of the form 
\begin{equation}\label{eq:examp.1}
u(x,y)=U\nabla^\perp\Psi_m(x,y)=UL_y\left(\frac{\partial\Psi_m}{\partial
    y},-\frac{\partial\Psi_m}{\partial x}\right)
\end{equation}
with the stream function 
\begin{equation}\label{eq:examp.2}
\Psi_m(x,y)=\cos^m(\pi x/L_x)\cos^m(\pi y/L_y),~~
m\ge 1
\end{equation}
periodic in $x$ and $y$. The scalar $U$ has the dimension of velocity.
The structure of the level sets of the functions $\Psi_m$ is clearly
the same for all $m$.  The period cell for these flows is the
rectangle $D=[-\frac{L_x}{2},\frac{3L_x}{2}] \times[-\frac{
  L_y}{2},\frac{3 L_y}{2}]$, that consists of four smaller rectangles
separated by separatrices $\Psi_1=0$. The normal component of $u(x,y)$
is equal to zero at the boundary of the period cell of $u(x,y)$, which
slows down the burning as compared to percolating flows. This effect
is quantified by the following Proposition.
\begin{prop}\label{thm:sublinear}
  Let $T(x,y,t)$ be the solution of the reaction diffusion equation
  (\ref{eq:1.1}) with either Neumann or periodic boundary conditions
  (\ref{eq:1.2}) and (\ref{eq:1.2.1}), respectively.  Let $u(x,y)$ be
  given by (\ref{eq:examp.1}), (\ref{eq:examp.2}) with $L_y=H/2.$
 Moreover, assume that the initial data $T_0(x,y)$ has the
  property that $T_0(x,y)=1$ for $x\le x_0$, and $T_0(x,y)=0$ for
  $x\ge x_1$.  Then there exists a constant $C>0$ such that
  for $U\ge v_0$ we have
\begin{equation}
  \label{eq:2.3.15}
  \frac{\langle V\rangle_\infty}{v_0} \le
C\left(1+\frac{l}{L_x}\right)\left(\frac{U}{v_0}\right)^{2/(1+m)}+
\frac{L_x}{4l}.
\end{equation} 
\end{prop} 
{\bf Proof.} We will construct a function $\phi(x)$, independent of
$y$ (and hence satisfying both the Neumann and periodic boundary
conditions (\ref{eq:1.2}) and (\ref{eq:1.2.1})), and $L_x$-periodic in
$x$, such that the function $\Phi(x,t)=e^{-\lambda(x-ct)}\phi(x)$
satisfies the inequality
\begin{equation}\label{eq:Phieq}
\Phi_t+u\cdot\nabla\Phi-\kappa\Delta\Phi-\frac{v_0^2}{4\kappa}\Phi\ge
0. 
\end{equation} 
Moreover, the function $\phi(x)$ will be positive, bounded, and bounded 
away from zero.
 Then maximum principle will imply
that the solution $T(x,y,t)$ of (\ref{eq:1.1}) with the Neumann or periodic
boundary conditions satisfies the inequality
\[ T(x,y,t)\le C_\lambda e^{-\lambda(x-ct)} 
\] 
since it holds at $t=0$ for all $\lambda >0$ because of our choice of the
initial data. Then we will have 
\begin{equation}\label{eq:vinftyc}
\langle V\rangle_\infty\le c 
\end{equation}
as in Lemma \ref{lemma0.1}. Therefore our goal is to find a function
$\phi(x)$ and $\lambda>0$ so as to satisfy (\ref{eq:Phieq}) with as
small $c$ as possible. The function $\phi(x)$ should obey the
inequality
\begin{eqnarray}
  \label{eq:2.2.4}
  L\phi=\frac{\kappa}{\lambda}\Delta\phi-2\kappa\phi_x-\frac{u}{\lambda}
\cdot\nabla\phi+u_1\phi\le B\phi
\end{eqnarray}
with
\begin{equation}\label{eq:2.2.4.1}
B=c-\kappa\lambda-\frac{v_0^2}{4\kappa\lambda}.
\end{equation}
We will define $\phi(x)$ on the interval $[-L_x/2,3L_x/2]$, and then
extend it periodically to the whole real line.  In order to make use
of the fact that the $x$-component of $u(x,y)$ is small near the
lines $x=-L_x/2,3L_x/2$ we consider a smooth cut-off function
$\chi(x)$ defined as follows. Let $\eta(x)$ be a cut-off function
\begin{eqnarray}
  \eta(x)=\left\{\begin{matrix}
{1, & |x|\le\frac{1}{2}\left(\frac{U}{v_0}\right)^{-\alpha} 
\cr
0, & 
|x|\ge \left(\frac{U}{v_0}\right)^{-\alpha},\cr}
\end{matrix}\right.\nonumber
\end{eqnarray}
and $\eta$ decays monotonically from one to zero between those
intervals, so that 
\[
|\eta'|\le C\left(\frac{U}{v_0}\right)^{\alpha},~~~~
|\eta''|\le C\left(\frac{U}{v_0}\right)^{2\alpha}.
\]
The exponent $\alpha >0$ is to be chosen later. Then we define for
$x\in[-L_x/2,3L_x/2]$
\[
\chi(x)=\eta\left(\frac{x}{L_x}+\frac{1}{2}\right)+
\eta\left(\frac{x}{L_x}-\frac{3}{2}\right),
\]
so that the two terms have non-overlapping support,
and set
\begin{eqnarray}
\phi(x)=\chi({x})+(1-\chi({x}))e^{\lambda x}:=\chi(x)+\beta(x)\nonumber
\end{eqnarray}
so that
\[
e^{-\lambda L_x/2}\le\phi(x,y)\le e^{3\lambda L_x/2}.
\]
First we observe that since $\displaystyle |u_1|\le
CU\left({U}/{v_0}\right)^{-m\alpha}$ on the support of $\chi(x)$, we have
\begin{eqnarray}
  \label{eq:2.2.8}
  |L\chi(x)|&=&\left|\frac{\kappa}{\lambda}\chi''-
2\kappa\chi'-\frac{u_1}{\lambda}\chi'+u_1\chi\right|\\&\le&
C\left[\frac{\kappa}{\lambda L_x^2}\left(\frac{U}{v_0}\right)^{2\alpha}+
\frac{\kappa}{L_x}\left(\frac{U}{v_0}\right)^\alpha+
\frac{v_0}{\lambda L_x}\left(\frac{U}{v_0}\right)^{1-m\alpha+\alpha}+
v_0\left(\frac{U}{v_0}\right)^{1-m\alpha}\right].
\nonumber
\end{eqnarray}
Moreover, we have
\begin{eqnarray}\label{eq:2.2.10}
|L\beta|=\left|-(1-\chi)\kappa\lambda
-\frac{\kappa}{\lambda}\chi''
+\frac{1}{\lambda}
u_1\chi'\right|e^{\lambda x}
\le C\left(
\frac{\kappa}{\lambda L_x^2}\left(\frac{U}{v_0}\right)^{2\alpha}
+\frac{v_0}{\lambda L_x}\left(\frac{U}{v_0}\right)^{1-m\alpha+\alpha}\right)
e^{\lambda x}+\kappa\lambda\beta.
\end{eqnarray}
We put together the bounds (\ref{eq:2.2.8}) and (\ref{eq:2.2.10}), and
obtain
\begin{eqnarray}
   \left|\frac{L\phi}{\phi}\right|\le Ce^{3\lambda L_x/2}
\left[
\frac{\kappa}{\lambda L_x^2}\left(\frac{U}{v_0}\right)^{2\alpha}+
\frac{\kappa}{L_x}\left(\frac{U}{v_0}\right)^\alpha+
v_0\left(\frac{U}{v_0}\right)^{1-m\alpha}+
\frac{v_0}{\lambda L_x}\left(\frac{U}{v_0}\right)^{1-m\alpha+\alpha}
\right]+\kappa\lambda. \label{eq:2.2.13}
\end{eqnarray}
Therefore the function $\phi$ that we have constructed satisfies the
inequality (\ref{eq:2.2.4}) with the constant $B$ given by the right
side of (\ref{eq:2.2.13}). Using the definition (\ref{eq:2.2.4.1}) of $B$
and relation (\ref{eq:vinftyc}) we obtain then the following bound on
the bulk burning rate:
\[
  \frac{\langle V\rangle_\infty}{v_0}\le 
Ce^{3\lambda L_x/2}\left[\frac{l}{\lambda
L_x^2}\left(\frac{U}{v_0}\right)^{2\alpha}+
\frac{l}{L_x}\left(\frac{U}{v_0}\right)^\alpha+
\left(\frac{U}{v_0}\right)^{1-m\alpha}+
\frac{1}{\lambda L_x}\left(\frac{U}{v_0}\right)^{1-m\alpha+\alpha}
\right]+
2l\lambda+\frac{1}{4l\lambda}
\]
where $l=\kappa/v_0$ is the laminar front width. 
Then we let $\lambda=1/L_x$ and $\alpha=1/(1+m)$,
and get the estimate in Proposition \ref{thm:sublinear}.$\Box$

One can see from the proof of Proposition \ref{thm:sublinear} that it
may be easily generalized to include cellular flows other than those of the
form (\ref{eq:examp.1}-\ref{eq:examp.2}). The relevant assumptions
are similar geometric structure of the streamlines and the appropriate
rate of decay of the normal velocity at the boundary of the period
cell. 

The power $\xi=2/(1+m)$ in (\ref{eq:2.3.15}) is probably not sharp,
but Proposition \ref{thm:sublinear} still shows several important
points. Not all advection velocities with non-trivial $u_1$ lead to linear growth of
the bulk burning rate in the advection amplitude. The presence of the 
closed streamlines appears to be crucial for sub-linear enhancement.
Also, the exponent $\xi$ may be
made arbitrarily close to $\xi=0$ by taking $m\to\infty$, and thus one
can construct non-trivial flows for which the bulk burning rate grows slower than
any given power of $U/v_0$.

\begin{appendix}

\section{Homogenization regime}\label{A1}

Here we briefly present a simple direct application of the bound we obtained
in Section \ref{glb}. We will consider a homogenization regime where 
the reaction is very weak, and investigate an effect of the periodic advection 
velocity in this limit. 
Let us consider the reaction-diffusion-advection equation
(\ref{eq:1.1}) with the laminar velocity $v_0$ being small:
$v_0\to\frac{1}{N}v_0$, $N\gg 1$. The domain is then taken to be finite
but very large: $D_N=ND$, where $D$ is some fixed region, such as a rectangle.
Initial data varies on the large scale:
\begin{eqnarray}
  \label{eq:2.10}
  &&(T_N)_t+u({\bf x},t)\cdot\nabla T_N =
\kappa\Delta T_N+\frac{v_0^2}{4N^2\kappa}f(T_N)\\
&&\frac{\partial T_N}{\partial n}=0~~\hbox{on}~~\partial D_N\nonumber\\
&&T_N({\bf x},0)=T_0(\frac{{\bf x}}{N}),~~{\bf x}\in D_N\nonumber
\end{eqnarray}
Then after rescaling ${\bf x},t\to N{\bf x},N^2 t$ the rescaled problem is
\begin{eqnarray}
  \label{eq:2.11}
  &&(T_N)_t+Nu(N{\bf x},N^2 t)\cdot\nabla T_N=
\kappa\Delta T_N+\frac{v_0^2}{4\kappa }f(T_N)\\
&&\frac{\partial T_N}{\partial n}=0~~\hbox{on}~~\partial D\nonumber\\
&&T_N({\bf x},0)=T_0({\bf x}),~~{\bf x}\in D_N\nonumber
\end{eqnarray}
We assume that $u({\bf x},t)$ is periodic in $x$ with period cell $Q$ and
in $t$ with period $\tau$, and
vanishes on the boundary of the cell $C$. Moreover, it is convenient to
assume that $D$ contains an integer number of cells, so that $u(N{\bf
  x},Nt)$ vanishes on the boundary $\partial D$. The bulk burning rate is
given as before by
\[
V_N(t)=\int\limits_{D}d{\bf x}\frac{\partial T_N}{\partial
  t}=\frac{v_0^2}{4\kappa}\int\limits_Dd{\bf x}f(T_N).
\]
The following Theorem may be established using the technique of \cite{BLP}.
\begin{thm}\label{thm3}
  The family of solutions $T_N$ of (\ref{eq:2.11}) converges strongly
  in $L_2([0,r]\times D)$ to the solution $\bar T$ of the homogenized
  problem
\begin{eqnarray}\label{eq:2.12}
&&\bar T_t=
\kappa_{ij}^*\frac{\partial^2 \bar T}{\partial x_i\partial x_j}+
\frac{v_0^2}{4\kappa}f(\bar T)\\
&&\frac{\partial \bar T}{\partial n}=0~~\hbox{on}~~\partial D\nonumber\\
&&\bar T({\bf x},0)=T_0({\bf x}).\nonumber
\end{eqnarray}
The tensor $\kappa^*$ is given by
\[
\kappa_{ij}^*=\kappa\delta_{ij}-\frac{1}{|Q|\tau}\int\limits_0^\tau ds\int\limits_Q
d{\bf y}u_i({\bf y},s)\theta_j({\bf y},s)
\]
with $\theta_i$ being the periodic solution of the cell problem
\[
\frac{\partial\theta_i}{\partial t}+u({\bf
  x},t)\cdot\nabla\theta_i-\kappa\Delta\theta_i =-u_i({\bf x},t).
\]
Moreover, there exists a constant $C$ such that
$\displaystyle\|T_N-\bar T\|_{L^2([0,r]\times D)}\le \frac{C}{N}$.
\end{thm}
Since $D$ is finite,
Theorem \ref{thm3} implies that the bulk burning rate $\displaystyle
V_N(t)\to \bar V(t)=\int_D\bar T_t(s)ds$.
Let us denote $k_*$ the minimal eigenvalue of the symmetric part
of the tensor $\kappa^*,$ and set $v_0^* = v_0 \sqrt{k_*/\kappa}.$  
This Theorem may be also
applied to the front propagation problem in a finite rectangle with 
the boundary conditions (\ref{eq:1.2}), (\ref{eq:1.3}), (\ref{eq:1.3a}). 
Arguments similar
to Theorem
\ref{thm2} imply that the bulk burning rate for the homogenized problem
obeys the lower bound
\[
\bar V(t)\ge Cv_0\sqrt{\frac{\beta k_*}{4\alpha\kappa}}(1-e^{-\alpha
v_0^2t/2\kappa})
\]
for times less than $t_* \approx \frac{\hbox{diam} D}{v^*_0}.$
Therefore the bulk burning rate in the original unscaled variables is
increased from $\frac{1}{N} v_0$ to
$\displaystyle\frac{1}{N}\sqrt{k_*/{\kappa}}$, but remains
of order $O(\frac{1}{N})$. The dependence of the tensor $\kappa^*$ on
the advection velocity $u$ and diffusivity $\kappa$ is rather
complicated.  Some estimates for $\kappa^*$ were obtained in \cite{AM}
and \cite{FP}, and they may be applied to obtain the relevant bounds
for $\bar V(t)$. This homogenization analysis is applicable only in
the limit of very weak reaction at a fixed diffusivity. This is the case
when the front width $l_0=\kappa/v_0$ is much larger than the typical scale of
variations of the turbulent velocity. 

\end{appendix}

\noindent {\bf Acknowledgments} This work was partially supported by 
the ASCI Flash Center at the University of Chicago under DOE contract 
B341495. PC acknowledges partial
support of NSF-DMS9802611. AK and LR acknowledge partial support of
 the University of Chicago NSF-MRSEC. 
The authors are grateful to Fausto Cattaneo, Andrea Malagoli, Takis Souganidis, 
Vladimir Volpert and Natalia Vladimirova for
valuable discussions.

\end{document}